\documentclass{amsart}

\DeclareMathSymbol{\twoheadrightarrow}  {\mathrel}{AMSa}{"10}

\def\GG{{\mathcal G}}
\def\HH{{\mathcal H}}
\def\H{{\mathbb H}}
\def\A{{\mathbb A}}
\def\Q{{\mathbb Q}}
\def\Z{{\mathbb Z}}
\def\C{{\mathbb C}}
\def\R{{\mathbb R}}
\def\F{{\mathbb F}}
\def\P{{\mathbb P}}

\def\Sn{{\mathbf S}_n}

\def\RR{{\mathfrak R}}
\def\MM{{\mathbf M}}
\def\Perm{\mathrm{Perm}}
\def\Gal{\mathrm{Gal}}
\def\PSL{\mathrm{PSL}}
\def\PGL{\mathrm{PGL}}
\def\End{\mathrm{End}}
\def\Aut{\mathrm{Aut}}
\def\Hom{\mathrm{Hom}}

\def\I{{\mathcal I}}
\def\ZZ{{\mathcal Z}}

\def\J{{\mathcal J}}

\def\fchar{\mathrm{char}}
\def\GL{\mathrm{GL}}
\def\SL{\mathrm{SL}}

\def\pr{\mathrm{pr}}
\def\M{\mathrm{M}}

\def\dim{\mathrm{dim}}
\def\Res{\mathrm{Res}}

\newtheorem{thm}{Theorem}[section]
\newtheorem{lem}[thm]{Lemma}
\newtheorem{cor}[thm]{Corollary}

\theoremstyle{definition}

\newtheorem{ex}[thm]{Example}
\newtheorem{rem}[thm]{Remark}

\title[Homomorphisms of abelian varieties]
{Homomorphisms of abelian varieties}
\author[Yuri\ G.\ Zarhin]{Yuri\ G.\ Zarhin}
\address{Department of Mathematics, Pennsylvania State University,
University Park, PA 16802, USA} \email{zarhin\char`\@math.psu.edu}

\begin{document}

\maketitle

It is well-known that an abelian variety  is (absolutely) simple
or is isogenous to a self-product of an (absolutely) simple
abelian variety if and only if the center of its endomorphism
algebra is a field. In this paper we prove that the center is a
field if the field of definition of points of prime order $\ell$
is ``big enough".

The paper is organized as follows. In \S \ref{endo} we discuss
Galois properties of points of  order $\ell$ on an abelian variety
$X$ that imply that its endomorphism algebra $\End^0(X)$ is a
central simple algebra over the field of rational numbers. In \S
\ref{tate} we prove that similar Galois properties for two abelian
varieties $X$ and $Y$ combined with the linear disjointness of the
corresponding fields of definitions of points of  order $\ell$
imply that $X$ and $Y$ are non-isogenous (and even $\Hom(X,Y)=0$).
In \S \ref{jac} we give applications to endomorphism algebras of
hyperelliptic jacobians.
 In \S \ref{mult} we prove that if $X$ admits multiplications by a
 number field $E$ and the dimension of the centralizer of $E$ in
 $\End^0(X)$ is ``as large as possible" then $X$ is an abelian
 variety of CM-type isogenous to a self-product of an absolutely
 simple abelian variety.

 Throughout the paper we will freely use the following
 observation
 \cite[p. 174]{MumfordAV}: if an abelian variety $X$
 is isogenous to a self-product $Z^d$ of an abelian variety $Z$
 then a choice of an isogeny between  $X$ and $Z^d$ defines an isomorphism
 between $\End^0(X)$ and the  algebra
 $\M_d(\End^0(Z))$ of $d\times d$ matrices over $\End^0(Z)$. Since
 the center of $\End^0(Z)$ coincides with the center of
 $\M_d(\End^0(Z))$,
 we get an isomorphism between the center of $\End^0(X)$ and  the center of $\End^0(Z)$
 (that does not depend on the choice of an isogeny). Also
 $\dim(X)=d\cdot\dim(Z)$; in particular, both $d$ and $\dim(Z)$ divide
 $\dim(X)$.

\section{Endomorphism algebras of abelian varieties}
\label{endo} Throughout this paper $K$ is a field. We write $K_a$
for its algebraic closure and $\Gal(K)$ for the absolute Galois
group $\Gal(K_a/K)$. We write $\ell$ for a prime different from
$\fchar(K)$. If $X$ is an abelian variety of positive dimension
 over $K_a$ then we write $\End(X)$ for the ring of all its
$K_a$-endomorphisms and $\End^0(X)$ for the corresponding
$\Q$-algebra $\End(X)\otimes\Q$. If $Y$ is (may be, another)
abelian variety over $K_a$ then we write $\Hom(X,Y)$ for the group
of all $K_a$-homomorphisms from $X$ to $Y$. It is well-known that
 $\Hom(X,Y)=0$ if and only if $\Hom(Y,X)=0$.

 If $n$ is a positive
integer that is not divisible by $\fchar(K)$ then we write $X_n$
for the kernel of multiplication by $n$ in $X(K_a)$. It is
well-known \cite{MumfordAV} that  $X_n$ is a free $\Z/n\Z$-module
of rank $2\dim(X)$. In particular, if $n=\ell$ is a prime then
$X_{\ell}$ is an $\F_{\ell}$-vector space of dimension $2\dim(X)$.

 If $X$ is defined over $K$ then $X_n$ is a Galois
submodule in $X(K_a)$. It is known that all points of $X_n$ are
defined over a finite separable extension of $K$. We write
$\bar{\rho}_{n,X,K}:\Gal(K)\to \Aut_{\Z/n\Z}(X_n)$ for the
corresponding homomorphism defining the structure of the Galois
module on $X_n$,
$$\tilde{G}_{n,X,K}\subset
\Aut_{\Z/n\Z}(X_{n})$$ for its image $\bar{\rho}_{n,X,K}(\Gal(K))$
and $K(X_n)$ for the field of definition of all points of $X_n$.
Clearly, $K(X_n)$ is a finite Galois extension of $K$ with Galois
group $\Gal(K(X_n)/K)=\tilde{G}_{n,X,K}$. If $n=\ell$  then we get
a natural faithful linear representation
$$\tilde{G}_{\ell,X,K}\subset \Aut_{\F_{\ell}}(X_{\ell})$$
of $\tilde{G}_{\ell,X,K}$ in the $\F_{\ell}$-vector space
$X_{\ell}$.

\begin{rem}
\label{ell2}
 If $n=\ell^2$ then there is the  natural surjective
homomorphism
$$\tau_{\ell,X}:\tilde{G}_{\ell^2,X,K}\twoheadrightarrow \tilde{G}_{\ell,X,K}$$
corresponding to the field inclusion $K(X_{\ell})\subset
K(X_{\ell^2})$; clearly, its kernel is a finite $\ell$-group.
Every prime dividing $\#(\tilde{G}_{\ell^2,X,K})$ either divides
$\#(\tilde{G}_{\ell,X,K})$ or is equal to $\ell$. If $A$ is a
subgroup in  $\tilde{G}_{\ell^2,X,K}$ of index $N$ then its image
$\tau_{\ell,X}(A)$ in $\tilde{G}_{\ell,X,K}$ is isomorphic to
$A/A\bigcap\ker(\tau_{\ell,X})$. It follows easily that the index
of  $\tau_{\ell,X}(A)$ in $\tilde{G}_{\ell,X,K}$ equals $N/\ell^j$
where $\ell^j$ is  the index of $A\bigcap\ker(\tau_{\ell,X})$ in
$\ker(\tau_{\ell,X})$. In particular, $j$ is a nonnegative
integer.
\end{rem}

We write $\End_K(X)$ for the ring of all $K$-endomorphisms of $X$.
We have
$$\Z=\Z\cdot 1_X \subset \End_K(X) \subset \End(X)$$ where $1_X$
is the identity automorphism of $X$. Since $X$ is defined over
$K$, one may associate with every $u \in \End(X)$ and $\sigma \in
\Gal(K)$ an endomorphism $^{\sigma}u\ \in \End(X)$ such that
$^{\sigma}u (x)=\sigma u(\sigma^{-1}x)$ for $x \in X(K_a)$ and we
get the  group homomorphism
$$\kappa_{X}: \Gal(K) \to \Aut(\End(X));
\quad\kappa_{X}(\sigma)(u)=\ ^{\sigma}u \quad \forall \sigma \in
\Gal(K),u \in \End(X).$$ It is well-known that $\End_K(X)$
coincides with the subring of $\Gal(K)$-invariants in $\End(X)$,
i.e., $\End_K(X)=\{u\in \End(X)\mid\  ^{\sigma}u\ =u \quad \forall
\sigma \in \Gal(K)\}$. It is also well-known that $\End(X)$
(viewed as a group with respect to addition) is a free commutative
group of finite rank and $\End_K(X)$ is its {\sl pure} subgroup,
i.e., the quotient $\End(X)/\End_K(X)$ is also  a free commutative
group of finite rank.
 All  endomorphisms of $X$ are defined over a finite separable
 extension of $K$.
More precisely \cite{Silverberg}, if $n\ge 3$ is a positive
integer not divisible by $\fchar(K)$ then all the endomorphisms of
$X$ are defined over $K(X_n)$; in particular,
$$\Gal(K(X_n)) \subset \ker(\kappa_{X})\subset \Gal(K).$$
 This implies that if
$\Gamma_K: =\kappa_{X}(\Gal(K)) \subset \Aut(\End(X))$ then there
exists a surjective homomorphism $\kappa_{X,n}:\tilde{G}_{n,X}
\twoheadrightarrow \Gamma_K$ such that the composition
$$\Gal(K)\twoheadrightarrow \Gal(K(X_{n})/K)=
\tilde{G}_{n,X}\stackrel{\kappa_{X,n}}{\twoheadrightarrow}
\Gamma_K$$ coincides with $\kappa_X$ and
$$\End_K(X)=\End(X)^{\Gamma_K}.$$
Clearly,  $\End(X)$ leaves invariant the subgroup $X_{\ell}\subset
X(K_a)$. It is well-known that $u\in \End(X)$ kills $X_{\ell}$
(i.e. $u(X_{\ell})=0$) if and only if $u \in \ell\cdot\End(X)$.
This gives us a natural embedding
$$\End_K(X)\otimes\Z/\ell\Z\subset \End(X)\otimes\Z/\ell\Z\
\hookrightarrow \End_{\F_{\ell}}(X_{\ell});$$  the image of
$\End_K(X)\otimes\Z/\ell\Z$ lies in the centralizer of the Galois
group, i.e., we get an embedding
$$\End_K(X)\otimes\Z/\ell\Z
\hookrightarrow
\End_{\Gal(K)}(X_{\ell})=\End_{\tilde{G}_{\ell,X,K}}(X_{\ell}).$$
The next easy assertion seems to be well-known (compare with Prop.
3 and its proof on pp. 107--108 in \cite{Mori2}) but quite useful.

\begin{lem}
\label{overK}
 If $\End_{\tilde{G}_{\ell,X,K}}(X_{\ell})=\F_{\ell}$
then $\End_K(X)=\Z$.
\end{lem}

\begin{proof}
It follows that the $\F_{\ell}$-dimension of
$\End_K(X)\otimes\Z/\ell\Z$ does not exceed $1$. This means that
the rank of the free commutative group $\End_K(X)$ does not exceed
$1$ and therefore is $1$. Since $\Z\cdot 1_X \subset \End_K(X)$,
it follows easily that $\End_K(X)=\Z\cdot 1_X=\Z$.
\end{proof}

\begin{lem}
\label{overK2}
 If $\End_{\tilde{G}_{\ell,X,K}}(X_{\ell})$ is a field
then $\End_K(X)$ has no zero divisors, i.e, $\End_K(X)\otimes\Q$
is a division algebra over $\Q$.
\end{lem}

\begin{proof}
It follows that $\End_K(X)\otimes\Z/\ell\Z$ is also a field and
therefore has no zero divisors. Suppose that $u,v$ are non-zero
elements of $\End_K(X)$ with $uv=0$. Dividing (if possible) $u$
and $v$ by suitable powers of $\ell$ in $\End_K(X)$, we may assume
that both $u$ and $v$ do not lie in $\ell \End_K(X)$ and induce
non-zero elements in $\End_K(X)\otimes\Z/\ell\Z$ with zero
product. Contradiction.
\end{proof}

 Let us put $\End^0(X):=\End(X)\otimes\Q$. Then
 $\End^0(X)$ is a semisimple
finite-dimensional $\Q$-algebra \cite[\S 21]{MumfordAV}.  Clearly,
the natural map $\Aut(\End(X)) \to \Aut(\End^0(X))$ is an
embedding. This allows us to view $\kappa_{X}$ as a homomorphism
$$\kappa_{X}: \Gal(K) \to \Aut(\End(X))\subset \Aut(\End^0(X)),$$
whose image coincides with $\Gamma_K\subset \Aut(\End(X))\subset
\Aut(\End^0(X))$;
  the subalgebra $\End^0(X)^{\Gamma_K}$  of
 $\Gamma_K$-invariants coincides with $\End_K(X)\otimes\Q$.

\begin{rem}
\label{split}
\begin{itemize}
\item[(i)]
 Let us split
the semisimple $\Q$-algebra $\End^0(X)$ into a finite direct
product $\End^0(X)= \prod_{s\in \I} D_s$
 of simple $\Q$-algebras $D_s$. (Here $\I$ is identified with the
 set of minimal two-sided ideals in $\End^0(X)$.)
Let $e_s$ be the identity element of $D_s$. One may view $e_s$ as
an idempotent in $\End^0(X)$. Clearly,
$$1_X=\sum_{s\in \I} e_s\in \End^0(X), \quad e_s e_t=0 \ \forall s\ne t.$$
There exists a positive integer $N$ such that all $N \cdot e_s$
lie in $\End(X)$. We write $X_s$ for the image $X_s:=(Ne_s) (X)$;
it is an abelian subvariety in $X$ of positive dimension.  The sum
map
$$\pi_X:\prod_s X_s \to X, \quad (x_s) \mapsto \sum_s x_s$$
is an isogeny. It is also clear that the intersection $D_s\bigcap
\End(X)$ leaves $X_s \subset X$ invariant. This gives us a natural
identification $D_s \cong \End^0(X_s)$.  One may easily check that
each $X_s$ is isogenous to a self-product of  (absolutely) simple
abelian variety and if $s\ne t$ then $\Hom(X_s,X_t)=0$.
\item[(ii)]
We write $C_s$ for the center of $D_s$. Then $C_s$ coincides with
the center of $\End^0(X_s)$ and is therefore either a totally real
number field of degree dividing $\dim(X_s)$ or a CM-field of
degree dividing $2\dim(X_s)$ \cite[p. 202]{MumfordAV}; the center
$C$ of $\End^0(X)$ coincides with $\prod_{s\in \I}C_s=\oplus_{s\in
S}C_s$.
\item[(iii)]
All the sets $$\{e_s\mid s\in \I\}\subset \oplus_{s\in \I}\Q\cdot
e_s\subset \oplus_{s\in \I}C_s=C$$ are stable under the Galois
action $\Gal(K)
\stackrel{\kappa_X}{\longrightarrow}\Aut(\End^0(X))$. In
particular, there is a continuous homomorphism from $\Gal(K)$ to
the group $\Perm(\I)$ of permutations of $\I$ such that its kernel
contains $\ker(\kappa_X)$ and
$$e_{\sigma(s)}=\kappa_X(\sigma)(e_s)=\ ^{\sigma}e_s ,
\ ^{\sigma}(C_s)\ =C_{\sigma(s)}, \ ^{\sigma}(D_s)=D_{\sigma(s)}
 \quad \forall \sigma\in \Gal(K), s\in \I.$$
It follows that $X_{\sigma(s)}=Ne_{\sigma(s)}(X)=\sigma (Ne_s
(X))=\sigma(X_s)$; in particular, abelian subvarieties $X_s$ and
$X_{\sigma(s)}$ have the same dimension and $u\mapsto\ ^{\sigma}u$
gives rise to an isomorphism of $\Q$-algebras
$\End^0(X_{\sigma(s)})\cong \End^0(X_s)$.
\item[(iv)]
If $J$ is a non-empty Galois-invariant subset in $\J$ then the sum
$\sum_{s\in J}Ne_s$ is Galois-invariant and therefore lies in
$\End_K(X)$. If $J'$ is another Galois-invariant subset of $\I$
that does not meet $J$ then $\sum_{s\in J}Ne_s$ also lies in
$\End_K(X)$ and $\sum_{s\in J}Ne_s \sum_{s\in J'}Ne_s=0$. Assume
that $\End_K(X)$ has no zero divisors. It follows that $\I$ must
consist of one Galois orbit; in particular, all $X_s$ have the
same dimension equal to $\dim(X)/\#(\I)$. In addition, if $t\in
\I$, $\Gal(K)_t$ is the stabilizer of $t$ in $\Gal(K)$ and $F_t$
is the subfield of $\Gal(K)_t$-invariants in the separable closure
of $K$ then it follows easily that $\Gal(K)_t$ is an open subgroup
of index $\#(\I)$ in $\Gal(K)$, the field extension $F_t/K$ is
separable  of degree $\#(\I)$ and $\prod_{s\in S}X_s$ is
isomorphic over $K_a$ to the Weil restriction $\Res_{F_t/K}(X_t)$.
This implies that $X$ is isogenous over $K_a$ to
$\Res_{F_t/K}(X_t)$.
\end{itemize}
\end{rem}

\begin{thm}
\label{center1} Suppose  that $\ell$ is a prime, $K$ is a field of
characteristic $\ne \ell$. Suppose that $X$ is an abelian variety
of positive dimension $g$ defined over  $K$.  Assume that
$\tilde{G}_{\ell,X,K}$ contains a  subgroup $\mathcal{G}$ such
 $\End_{\mathcal{G}}(X_{\ell})$ is a field.

Then one of the following  conditions holds:

\begin{enumerate}
\item[(a)]
The center of  $\End^0(X)$ is a field. In other
 words, $\End^0(X)$ is a simple $\Q$-algebra.
\item[(b)]
\begin{itemize}
\item[(i)]
The prime $\ell$ is odd;
\item[(ii)]
there exist a positive integer $r>1$ dividing $g$, a field $F$
with
$$K\subset K(X_{\ell})^{\mathcal{G}}=:L\subset F\subset K(X_{\ell}), \quad [F:L]=r$$
 and a $\frac{g}{r}$-dimensional abelian variety
$Y$ over $F$ such that $\End^0(Y)$ is a simple $\Q$-algebra, the
$\Q$-algebra $\End^0(X)$ is isomorphic to the direct sum of $r$
copies of $\End^0(Y)$ and the Weil restriction $\Res_{F/L}(Y)$ is
isogenous over $K_a$ to $X$. In particular, $X$ is isogenous over
$K_a$ to a product of $\frac{g}{r}$-dimensional abelian varieties.
In addition, $\mathcal{G}$ contains a subgroup of index $r$;
\end{itemize}
 \item[(c)]
\begin{itemize}
\item[(i)]
The prime $\ell=2$;
\item[(ii)]
there exist a positive integer $r>1$ dividing $g$,  fields $L$ and
$F$ with
$$K\subset K(X_{4})^{\mathcal{G}}\subset L\subset F\subset K(X_{4}), \quad [F:L]=r$$ and a
$\frac{g}{r}$-dimensional abelian variety $Y$ over $F$ such that
$\End^0(Y)$ is a simple $\Q$-algebra,  the $\Q$-algebra
$\End^0(X)$ is isomorphic to the direct sum of $r$ copies of
$\End^0(Y)$ and the Weil restriction $\Res_{F/L}(Y)$ is isogenous
over $K_a$ to $X$. In particular, $X$ is isogenous over $K_a$ to a
product of $\frac{g}{r}$-dimensional abelian varieties.In
addition, there exists a nonnegative integer $j$ such that $2^j$
divides $r$ and $\mathcal{G}$ contains a subgroup of index
$\frac{r}{2^j}>1$.
\end{itemize}
\end{enumerate}
\end{thm}

\begin{proof}
We will use notations of Remark \ref{split}.  Let us put $n=\ell$
if $\ell$ is odd and $n=4$ if $\ell=2$. Replacing $K$ by
$K(X_{\ell})^{{\mathcal{G}}}$, we may and will assume that
$$\tilde{G}_{\ell,X,K}={\mathcal{G}}.$$
If $\ell$ is odd then let us put $L=K$ and
$H:=\Gal(K(X_{\ell})/K)={\mathcal{G}}=\Gal(L(X_{\ell})/L)$.

If $\ell=2$ then we choose a subgroup $\HH \subset
\tilde{G}_{4,X,K}$ of smallest possible order such that
$\tau_{2,X}(\HH)=\tilde{G}_{2,X,K}={\mathcal{G}}$ and put
$L:=K(X_{4})^{\HH}\subset K(X_{4})$. It follows easily that
$L(X_{4})=K(X_{4})$ and
 $\Gal(L(X_{2})/L)=\Gal(K(X_{2})/K)$, i.e.,
$$\HH=\tilde{G}_{4,X,L},\quad
\tilde{G}_{2,X,L}={\mathcal{G}}.$$ The minimality property of
$\HH$ combined with Remark \ref{ell2} implies that if $H\subset
\tilde{G}_{4,X,L}$ is a subgroup of index $r>1$ then
$\tau_{2,X}(H)$ has index $\frac{r}{2^j}>1$ in $\tilde{G}_{2,X,L}$
for some nonnegative index $j$.

In light of Lemma \ref{overK2}, $\End_L(X)$ has no zero divisors.
It follows from  Remark \ref{split}(iv) that $\Gal(L)$ acts on
$\I$ transitively. Let us put $r=\#(\I)$. If $r=1$ then $\I$ is a
singleton and $\I=\{s\}, X=X_s, \End^0(X)=D_s,C=C_s$. This means
that  assertion (a) of Theorem \ref{center1} holds true.

Further we assume that $r>1$. Let  us choose $t\in \I$ and put
$Y:=X_t$. If $F:=F_t$ is the subfield of $\Gal(L)_t$-invariants in
the separable closure of $K$ then it follows from  Remark
\ref{split}(iv) that $F_t/L$ is a separable degree $r$ extension,
$Y$ is defined over $F$ and $X$ is isogenous over $L_a=K_a$ to
$\Res_{F/L}(Y)$.

Recall (Remark \ref{split}(iii)) that $\ker(\kappa_X)$ acts
trivially on $\I$. It follows that  $\Gal(L(X_n))$ acts trivially
on $\I$. This implies that $\Gal(L(X_n))$ lies in $\Gal(L)_t$.
Recall that $\Gal(L)_t$ is an open subgroup of index $r$ in
$\Gal(L)$ and $\Gal(L(X_n))$ is a normal open subgroup in
$\Gal(L)$. It follows that $H:=\Gal(L)_t/\Gal(L(X_n))$ is a
subgroup of index $r$ in
$$\Gal(L)/\Gal(L(X_n))=\Gal(L(X_n)/L)=\tilde{G}_{n,X,L}.$$
If $\ell$ is odd then $n=\ell$ and
$\tilde{G}_{n,X,L}=\tilde{G}_{\ell,X,L}=\GG$  contains a subgroup
of index $r>1$. It follows from Remark \ref{split} that assertion
(b) of Theorem \ref{center1} holds true.

If $\ell=2$ then $n=4$ and $\tilde{G}_{n,X,L}=\tilde{G}_{4,X,L}$
contains a subgroup $H$ of index $r>1$. But in this case we know
(see the very beginning of this proof) that
$\tilde{G}_{2,X,L}=\GG$ and $\tau_{2,X}(H)$ has index
$\frac{r}{2^j}>1$ in $\tilde{G}_{2,X,L}$ for some nonnegative
integer $j$. It follows from Remark \ref{split} that  assertion
(c) of Theorem \ref{center1} holds true.
\end{proof}

Before stating our next result, recall that a {\sl perfect} finite
group $\GG$ with center $\ZZ$ is called {\sl quasi-simple} if the
quotient $\GG/\ZZ$ is a simple nonabelian group. Let $H$ be a
non-central normal subgroup in quasi-simple ${\mathcal{G}}$. Then
the image of $H$ in simple ${\mathcal{G}}/\ZZ$ is a non-trivial
normal subgroup and therefore coincides with ${\mathcal{G}}/\ZZ$.
This means that ${\mathcal{G}}=\ZZ H$. Since ${\mathcal{G}}$ is
perfect, ${\mathcal{G}}=[{\mathcal{G}},{\mathcal{G}}]=[H,H]\subset
H$. It follows that ${\mathcal{G}}=H$. In other words,  every
proper normal subgroup in a quasi-simple group is central.

\begin{thm}
\label{center2} Suppose  that $\ell$ is a prime, $K$ is a field of
characteristic different from $\ell$. Suppose that $X$ is an
abelian variety of positive dimension $g$ defined over  $K$.   Let
us assume that $\tilde{G}_{\ell,X,K}$ contains a subgroup
$\mathcal{G}$ that enjoys the following properties:

\begin{enumerate}
\item[(i)]
$\End_{\mathcal{G}}(X_{\ell})=\F_{\ell}$;
\item[(ii)]
The group $\mathcal{G}$ does not contain a subgroup of index $2$.
\item[(iii)]
The only normal subgroup in $\mathcal{G}$ of index dividing $g$ is
$\mathcal{G}$ itself.
\end{enumerate}

 Then one of the following two conditions (a) and (b) holds:
\begin{enumerate}
\item[(a)] There exists a positive integer $r>2$ such that:
\begin{itemize}
\item[(a0)]
$r$ divides $g$ and $X$ is isogenous over $K_a$ to a product of
$\frac{g}{r}$-dimensional abelian varieties;
\item[(a1)]
If $\ell$ is odd then
 $\mathcal{G}$  contains a subgroup of index $r$;
\item[(a2)]
If $\ell=2$ then there exists a nonnegative integer $j$ such that
 $\mathcal{G}$  contains a subgroup of index $\frac{r}{2^j}>1$.
\end{itemize}

\item[(b)]
\begin{itemize}
\item[(b1)]
The center of  $\End^0(X)$ coincides with $\Q$. In other
 words, $\End^0(X)$ is a matrix algebra either over $\Q$ or over a
 quaternion $\Q$-algebra.
\item[(b2)]
 If $\mathcal{G}$ is perfect and $\End^0(X)$ is a matrix algebra  over
 a quaternion $\Q$-algebra $\H$ then $\H$ is unramified at every
 prime not dividing  $\#(\mathcal{G})$.
\item[(b3)]
Let $\ZZ$ be the center of $\mathcal{G}$.  Suppose that
$\mathcal{G}$ is quasi-simple, i.e. it is perfect and  the
quotient ${\mathcal{G}}/\ZZ$ is a simple group. If $\End^0(X) \ne
\Q$ then there exist a perfect finite (multiplicative) subgroup
$\Pi\subset\End^0(X)^{*}$ and a surjective homomorphism $\Pi
\twoheadrightarrow {\mathcal{G}}/\ZZ$.
 \end{itemize}
 \end{enumerate}
\end{thm}

\begin{proof}
Let $C$ be the center of  $\End^0(X)$. Assume that  $C$  is not a
field. Applying Theorem \ref{center1}, we conclude that the
condition (a) holds.

 Assume now that  $C$   is  a field.
 We need to prove (b).
Let us define $n$ and $L$ as in the beginning of the proof of
Theorem \ref{center1}. We have
$$\GG=\tilde{G}_{\ell,X,L},\quad
\End_{\tilde{G}_{\ell,X,L}}(X_{\ell})=\F_{\ell}.$$ In addition, if
$\ell=2$ and  $H\subset \tilde{G}_{4,X,L}$ is a subgroup of index
$r>1$ then $\tau_{2,X}(H)$ has index $\frac{r}{2^j}>1$ in
$\tilde{G}_{2,X,L}=\GG$ for some nonnegative integer $j$. This
implies that the only normal subgroup in
$\tilde{G}_{n,X,L}=\tilde{G}_{4,X,L}$ of index dividing $g$ is
$\tilde{G}_{n,X,L}$ itself. It is also clear that
$\tilde{G}_{n,X,L}$ does not contain a subgroup of index $2$. It
follows from Remark \ref{ell2} that if $\GG$ is perfect then
$\tilde{G}_{4,X,L}$ is also perfect and every prime dividing
$\#(\tilde{G}_{4,X,L})$ must divide $\#(\GG)$, because (thanks to
a celebrated theorem of Feit-Thompson) $\#(\GG)$ must be even. (If
$\ell$ is odd then $n=\ell$ and $\tilde{G}_{n,X,L}=\GG$.)

It follows from Lemma \ref{overK} that $\End_L(X)=\Z$ and
therefore $\End_L(X)\otimes\Q=\Q$. Recall that
$\End_L(X)\otimes\Q=\End^0(X)^{\Gal(L)}$ and
$\kappa_X:\Gal(L)\to\Aut(\End^0(X))$ kills $\Gal(L(X_n))$. This
gives rise to the homomorphism
$$\kappa_{X,n}:\tilde{G}_{n,X,L}=\Gal(L(X_n)/L)=\Gal(L)/\Gal(L(X_n))\to\Aut(\End^0(X))$$
with
$\kappa_{X,n}(\tilde{G}_{n,X,L})=\kappa_X(\Gal(L))\subset\Aut(\End^0(X))$
and $\End^0(X)^{\tilde{G}_{n,X,L}}=\Q$. Clearly, the action of
$\tilde{G}_{n,X,L}$ on $\End^0(X)$ leaves invariant the center $C$
and therefore defines a homomorphism $\tilde{G}_{n,X,L}\to
\Aut(C)$ with $C^{\tilde{G}_{n,X,L}}=\Q$. It follows that $C/\Q$
is a Galois extension and the corresponding map
$$\tilde{G}_{n,X,L}\to \Aut(C)=\Gal(C/\Q)$$
is surjective.
 Recall  that $C$ is either a totally real
number field of degree dividing $g$ or a purely imaginary
quadratic extension of a totally real number field $C^{+}$ where
$[C^{+}:\Q]$ divides $g$ . In the case of totally real $C$ let us
put $C^{+}:=C$. Clearly, in both cases $C^{+}$ is the largest
totally real subfield of $C$ and therefore the action of
$\tilde{G}_{n,X,L}$ leaves $C^{+}$ stable, i.e. $C^{+}/\Q$ is also
a Galois extension.  Let us put $r:=[C^{+}:\Q]$. It is known
\cite[ p. 202]{MumfordAV}
 that $r$ divides $g$.  Clearly, the Galois
group $\Gal(C^{+}/\Q)$ has order $r$ and we have a surjective
homomorphism (composition)
$$\tilde{G}_{n,X,L}\twoheadrightarrow
\Gal(C/\Q)\twoheadrightarrow\Gal(C^{+}/\Q)$$ of
$\tilde{G}_{n,X,L}$ onto order $r$ group $\Gal(C^{+}/\Q)$.
Clearly, its kernel is a normal subgroup of index $r$ in
$\tilde{G}_{n,X,L}$. This contradicts our assumption if $r>1$.
Hence $r=1$, i.e. $C^{+}=\Q$.
  It follows that either $C=\Q$ or $C$ is an
imaginary quadratic field and $\Gal(C/\Q)$ is a group of order
$2$. In the latter case we get the surjective homomorphism from
$\tilde{G}_{n,X,L}$ onto $\Gal(C/\Q)$, whose kernel is a subgroup
of order 2 in $\tilde{G}_{n,X,L}$, which does not exist. This
proves that $C=\Q$. It follows from Albert's classification
\cite[p. 202]{MumfordAV} that $\End^0(X)$ is either a matrix
algebra $\Q$ or a matrix algebra $\M_d(\H)$ where $\H$ is a
quaternion $\Q$-algebra. This proves  assertion (b1) of Theorem
\ref{center2}.

Assume, in addition, that $\GG$ is perfect. Then, as we have
already seen, $\tilde{G}_{n,X,L}$ is also perfect. This implies
that $\Gamma:=\kappa_{X,n}(\tilde{G}_{n,X,L})$ is a finite perfect
subgroup of $\Aut(\End^0(X))$ and every prime dividing
$\#(\Gamma)$ must divide $\#(\tilde{G}_{n,X,L})$ and therefore
divides $\#(\GG)$. Clearly,
$$\Q=\End^0(X)^{\Gamma} \eqno(1).$$
Assume that $\End^0(X)\ne\Q$. Then $\Gamma\ne \{1\}$. Since
$\End^0(X)$ is a  central simple $\Q$-algebra, all its
automorphisms are inner, i.e.,
$\Aut(\End^0(X))=\End^0(X)^*/\Q^{*}$. Let $\Delta
\twoheadrightarrow \Gamma$ be the universal central extension of
$\Gamma$. It is well-known \cite[Ch. 2, \S 9]{Suzuki} that
$\Delta$ is a finite perfect group and the set of prime divisors
of $\#(\Delta)$ coincides with the set of prime divisors of
$\#(\Gamma)$ . The universality property implies that the
inclusion map $\Gamma\subset \End^0(X)^{*}/\Q^{*}$ lifts
(uniquely) to a homomorphism $\pi:\Delta \to \End^0(X)^{*}$. The
equality (1) means that the centralizer of $\pi(\Delta)$ in
$\End^0(X)$ coincides with $\Q$ and therefore $\ker(\pi)$ does not
coincide with $\Delta$. It follows that the image $\Gamma_0$ of
$\ker(\pi)$ in $\Gamma$ does not coincide with the whole $\Gamma$.
It also follows that if $\Q[\Delta]$ is the group $\Q$-algebra of
$\Delta$ then $\pi$ induces the $\Q$-algebra homomorphism $\pi:
\Q[\Delta] \to \End^0(X)$ such that the centralizer of the image
$\pi(\Q[\Delta])$ in $\End^0(X)$ coincides with $\Q$.

I claim that $\pi(\Q[\Delta])=\End^0(X)$ and therefore $\End^0(X)$
is isomorphic to a direct summand of $\Q[\Delta]$. This claim
follows easily from the next lemma that will be proven later in
this section.

\begin{lem}
\label{sur} Let $E$ be a field of characteristic zero, $T$ a
semisimple finite-dimensional $E$-algebra, $S$ a
finite-dimensional central simple $E$-algebra, $\beta: T \to S$ an
$E$-algebra homomorphism that sends $1$ to $1$. Suppose that the
centralizer of the image $\beta(T)$ in $S$ coincides with the
center $E$. Then $\beta$ is surjective, i. e. $\beta(T)=S$.
\end{lem}

In order to prove  (b2), let us assume that $\End^0(X)=\M_d(\H)$
where $\H$ is a quaternion $\Q$-algebra. Then
 $\M_d(\H)$ is isomorphic to a direct summand of $\Q[\Delta]$. On
the other hand, it is well-known that if $q$ is a prime not
dividing  $\#(\Delta)$ then
$\Q_q[\Delta]=\Q[\Delta]\otimes_{\Q}\Q_q$ is a direct sum of
matrix algebras over (commutative) fields. It follows that
$\M_d(\H)\otimes_{\Q}\Q_q$  also splits. This proves the assertion
(b2).

In order to prove  (b3), let us assume that ${\mathcal{G}}$ is a
quasi-simple finite group with center $\ZZ$. Let us put
$\Pi:=\pi(\Delta)\subset \End^0(X)^{*}$. We are going to construct
a surjective homomorphism $\Pi \twoheadrightarrow
{\mathcal{G}}/\ZZ$. In order to do that, it suffices to construct
a surjective homomorphism $\Gamma \twoheadrightarrow
{\mathcal{G}}/\ZZ$. Recall that there are surjective homomorphisms
$$\tau:\tilde{G}_{n,X,L} \twoheadrightarrow
\tilde{G}_{\ell,X,L}={\mathcal{G}},\quad
\kappa_{X,n}:\tilde{G}_{n,X,L} \twoheadrightarrow\Gamma.$$ (If
$\ell$ is odd then $\tau$ is the identity map; if $\ell=2$ then
$\tau=\tau_{2,X}$.) Let $H_0$ be the kernel of
$\kappa_{X,n}:\tilde{G}_{n,X,L} \twoheadrightarrow\Gamma$.
Clearly,
$$\tilde{G}_{n,X,L}/H_0\cong \Gamma \eqno(2).$$ Since $\Gamma\ne \{1\}$, we have $H_0\ne
\tilde{G}_{n,X,L}$. It follows  that $\tau(H_0)\ne {\mathcal{G}}$.
The surjectivity of $\tau:\tilde{G}_{n,X,L} \twoheadrightarrow
{\mathcal{G}}$ implies that $\tau(H_0)$ is normal in
${\mathcal{G}}$ and therefore lies in the center $\ZZ$. This gives
us the surjective homomorphisms
$$\tilde{G}_{n,X,L}/H_0 \twoheadrightarrow
\tau(\tilde{G}_{n,X,L})/\tau(H_0)={\mathcal{G}}/\tau(H_0)\twoheadrightarrow
{\mathcal{G}}/\ZZ,$$ whose composition is a surjective
homomorphism $\tilde{G}_{n,X,L}/H_0\twoheadrightarrow
{\mathcal{G}}/\ZZ$. Using (2), we get the desired surjective
homomorphism $\Gamma\twoheadrightarrow {\mathcal{G}}/\ZZ$.
\end{proof}

\begin{proof}[Proof of Lemma \ref{sur}]
Replacing $E$ by its algebraic closure $E_a$ and tensoring $T$ and
$S$ by $E_a$, we may and will assume that $E$ is algebraically
closed. Then $S=\M_n(E)$ for some positive integer $n$. Clearly,
$\beta(T)$ is a direct sum of say, $b$ matrix algebras over $E$
and the center of $\beta(T)$ is isomorphic to a direct sum of $b$
copies of $E$. In particular, if $b>1$ then the centralizer of
$\beta(T)$ in $S$ contains the $b$-dimensional center of
$\beta(T)$ which gives us the contradiction. So, $b=1$ and
$\beta(T)\cong \M_k(E)$ for some positive integer $k$. Clearly, $k
\le n$; if the equality holds then we are done. Assume that $k<n$:
we need to get a contradiction. So, we have
$$1 \in E \subset \beta(T)\cong \M_k(E)\hookrightarrow \M_n(E)=S.$$
This provides $E^n$ with a structure of faithful $\beta(T)$-module
in such a way that $E^n$ does not contain a non-zero submodule
with trivial (zero) action of $\beta(T)$. Since $\beta(T)\cong
\M_k(E)$, the $\beta(T)$-module $E^n$ splits into a direct sum of
say, $e$  copies of a simple faithful $\beta(T)$-module $W$ with
$\dim_E(W)=k$. Clearly, $e=n/k>1$. It follows easily that the
centralizer of $\beta(T)$ in $S=\M_n(E)$ coincides with
$$\End_{\beta(T)}(W^e)=\M_e(\End_{\beta(T)}(W))=\M_e(E)$$
and has $E$-dimension $e^2>1$. Contradiction.
\end{proof}

\begin{cor}
\label{odd} Suppose  that $\ell$ is a prime, $K$ is a field of
characteristic different from $\ell$. Suppose that $X$ is an
abelian variety of positive dimension $g$ defined over  $K$.   Let
us assume that $\tilde{G}_{\ell,X,K}$ contains a perfect subgroup
$\mathcal{G}$ that enjoys the following properties:
\begin{enumerate}
\item[(a)]
$\End_{\mathcal{G}}(X_{\ell})=\F_{\ell}$;
\item[(b)]
The only  subgroup of index dividing $g$ in $\mathcal{G}$ is
$\mathcal{G}$ itself.
\end{enumerate}
If $g$ is odd then either $\End^0(X)$ is a matrix algebra over
$\Q$ or $p=\fchar(K)>0$ and $\End^0(X)$ is a matrix algebra
$\M_d(\H_p)$
 over a quaternion $\Q$-algebra $\H_p$  that is
ramified exactly at $p$ and $\infty$ and $d>1$. In particular, if
$\fchar(K)$ does not divide $\#(\mathcal{G})$ then $\End^0(X)$ is
a matrix algebra over $\Q$.
\end{cor}
\begin{proof}[Proof of Corollary \ref{odd}]
Let us assume that $\End^0(X)$ is {\sl not} isomorphic to a matrix
algebra over $\Q$. Then $\End^0(X)$ is (isomorphic to) a matrix
algebra $\M_d(\H)$ over
 a quaternion $\Q$-algebra $\H$. This means that there exists an
 absolutely simple abelian variety $Y$ over $K_a$ such that $X$ is
 isogenous to $Y^d$ and $\End^0(Y)=\H$. Clearly, $\dim(Y)$ is odd.
It follows from Albert's classification \cite[p. 202]{MumfordAV}
that $p:=\fchar(K_a)=\fchar(K)> 0$.  By Lemma 4.3 of \cite{Oort},
if there exists a prime $q\ne p$ such that $\H$ is unramified at
$q$ then $4=\dim_{\Q}\H$ divides $2\dim(Y)$. Since $\dim(Y)$ is
odd, $2\dim(Y)$ is not divisible by $4$ and therefore $\H$ is
unramified at all primes different from $p$. It follows from the
theorem of Hasse-Brauer-Noether  that $\H \cong \H_p$.

Now, assume that $d=1$, i.e. $\End^0(X)=\H_p$. We know that
$\End^0(X)^{*}=\H_p^*$ contains a nontrivial finite perfect group
$\Pi$. But this contradicts to the following elementary statement,
whose proof will be given later in this section.

\begin{lem}
\label{Hp}
 Every finite subgroup in $\H_p^*$ is
solvable.
\end{lem}

Hence $\End^0(X)\ne \H_p$, i.e. $d>1$.

 Assume now that $p$ does {\sl not} divide
$\#(\mathcal{G})$. It follows from Theorem \ref{center2} that $\H$
is unramified at $p$. This implies that $\H$ can be ramified only
at $\infty$ which could not be the case. The obtained
contradiction proves that $\End^0(X)$ is a matrix algebra over
$\Q$.
\end{proof}

\begin{proof}[Proof of Lemma \ref{Hp}]
If $p\ne 2$ then $\H_p^*\subset (\H_p\otimes_{\Q}\Q_2)^{*}\cong
\GL(2,\Q_2)$ and if $p=2$ then $\H_2^*\subset
(\H_2\otimes_{\Q}\Q_3)^{*}\cong \GL(2,\Q_3)$. Since every finite
subgroup in $\GL(2,\Q_2)$ (resp. $\GL(2,\Q_3)$) is conjugate to a
finite subgroup in $\GL(2,\Z_2)$ (resp. $\GL(2,\Z_3)$), it
suffices to check that every finite subgroup in $\GL(2,\Z_2)$ and
$\GL(2,\Z_3)$ is solvable.

Recall that both $\GL(2,\F_2)$ and $\GL(2,\F_3)$ are solvable and
use the Minkowski-Serre lemma (\cite[pp. 124--125]{SerreL}; see
also \cite{SZ}).  This lemma asserts, in particular, that if $q$
is an odd prime  then the kernel of the reduction map
$\GL(n,\Z_q)\to\GL(n,\F_q)$ does not contain nontrivial elements
of finite order and that all periodic elements in the kernel of
the reduction map $\GL(n,\Z_2)\to\GL(n,\F_2)$ have order $1$ or
$2$.

 Indeed, every
finite subgroup $\Pi\subset \GL(2,\Z_3)$ maps injectively in
$\GL(2,\F_3)$ and therefore is solvable. If $\Pi\subset
\GL(2,\Z_2)$ is a finite subgroup then the kernel of the reduction
map $\Pi\to \GL(2,\F_2)$ consists of elements of order $1$ or $2$
and therefore is an elementary commutative $2$-group. Since the
image of the reduction map is solvable, we conclude that $\Pi$ is
solvable.
\end{proof}

\begin{cor}
\label{int} Suppose  that $\ell$ is a prime, $K$ is a field of
characteristic different from $\ell$. Suppose that $X$ is an
abelian variety of dimension $g$ defined over  $K$.  Let us put
$g'=\max(2,g)$. Let us assume that $\tilde{G}_{\ell,X,K}$ contains
a perfect subgroup $\mathcal{G}$ that enjoys the following
properties:
\begin{enumerate}
\item[(a)]
$\End_{\mathcal{G}}(X_{\ell})=\F_{\ell}$;
\item[(b)]
The only  subgroup of index dividing $g$ in $\mathcal{G}$ is
$\mathcal{G}$ itself.
\item[(c)] If  $\ZZ$ is the center of $\mathcal{G}$
 then ${\mathcal{G}}/\ZZ$ is a simple nonabelian group.
\end{enumerate}
Suppose that $\End^0(X)\cong \M_d(\Q)$ with $d>1$. Then there
exist a perfect finite  subgroup $\Pi\subset \GL(d,\Z)$ and a
surjective homomorphism $\Pi \twoheadrightarrow
{\mathcal{G}}/\ZZ$.
\end{cor}
\begin{proof}[Proof of Corollary \ref{int}]
Clearly, $\End^0(X)^{*}=\GL(n,\Q)$. One has only to recall that
every finite subgroup in $\GL(n,\Q)$ is conjugate to a finite
subgroup in $\GL(n,\Z)$ \cite[p. 124]{SerreL} and apply Theorem
\ref{center2}(iii).
\end{proof}

\section{Homomorphisms of abelian varieties}

\label{tate}
\begin{thm}
 \label{hom}
 Let $\ell$ be a prime, $K$ a field of characteristic different from
 $\ell$, $X$ and $Y$   abelian varieties of positive dimension defined over $K$.
 Suppose that the following conditions hold:

\begin{enumerate}
\item[(i)]
The extensions $K(X_{\ell})$ and $K(Y_{\ell})$  are linearly
disjoint over $K$.
\item[(ii)]
 $\End_{\tilde{G}_{\ell,X,K}}(X_{\ell})=\F_{\ell}$.
\item[(iii)]
The centralizer of $\tilde{G}_{\ell,Y,K}$ in
$\End_{\F_{\ell}}(Y_{\ell})$ is a field.
\end{enumerate}

Then either $\Hom(X,Y)=0,  \Hom(Y,X)=0$ or $\fchar(K)>0$ and both
abelian varieties $X$ and $Y$ are supersingular.
\end{thm}

\begin{rem} Theorem \ref{hom} was proven in \cite{ZarhinSh}  under  an
addititional assumption that the Galois modules $X_{\ell}$ and
$Y_{\ell}$ are simple. \end{rem}

In order to prove Theorem \ref{hom},  we need first to discuss the
notion of Tate module. Recall \cite{MumfordAV,SerreAb,ZP} that
this is a $\Z_{\ell}$-module $T_{\ell}(X)$
 defined as
the projective limit of Galois modules $X_{\ell^m}$. It is
well-known that $T_{\ell}(X)$ is a free $\Z_{\ell}$-module of rank
$2\dim(X)$ provided with the continuous action
$$\rho_{\ell,X}:\Gal(K) \to \Aut_{\Z_{\ell}}(T_{\ell}(X)).$$
There is  the natural isomorphism of
Galois modules
$$X_{\ell}=T_{\ell}(X)/\ell T_{\ell}(X) \eqno(3),$$
so one may view $\tilde{\rho}_{\ell,X}$ as the reduction of
$\rho_{\ell,X}$ modulo $\ell$.
 Let us put
$$V_{\ell}(X)=T_{\ell}(X)\otimes_{\Z_{\ell}}\Q_{\ell};$$
it is a $2\dim(X)$-dimensional $\Q_{\ell}$-vector space. The group
$T_{\ell}(X)$ is naturally identified with the $\Z_{\ell}$-lattice
in $V_{\ell}(X)$
 and the inclusion
$\Aut_{\Z_{\ell}}(T_{\ell}(X))\subset
\Aut_{\Q_{\ell}}(V_{\ell}(X))$ allows us to view $V_{\ell}(X)$ as
representation of $\Gal(K)$ over $\Q_{\ell}$. Let $Y$ be (may be,
another) abelian variety of positive dimension defined over $K$.
Recall \cite[\S 19]{MumfordAV}  that $\Hom(X,Y)$ is a free
commutative group of finite rank. Since $X$ and $Y$
 are defined over $K$, one may associate with every $u \in \Hom(X,Y)$
and $\sigma \in \Gal(K)$ an endomorphism $^{\sigma}u\ \in
\Hom(X,Y)$ such that
$$^{\sigma}u (x)=\sigma u(\sigma^{-1}x)\quad \forall x \in X(K_a)$$
and we get the group homomorphism
$$\kappa_{X,Y}: \Gal(K) \to \Aut(\Hom(X,Y));
\quad\kappa_{X,Y}(\sigma)(u)=\ ^{\sigma}u \quad \forall \sigma \in
\Gal(K),u \in \Hom(X,Y),$$ which provides the finite-dimensional
$\Q_{\ell}$-vector space $\Hom(X,Y)\otimes\Q_{\ell}$ with the
natural structure of Galois module.

There is a natural structure of  Galois module on the
$\Q_{\ell}$-vector space
$\Hom_{\Q_{\ell}}(V_{\ell}(X),V_{\ell}(Y))$  induced by the Galois
actions on  $V_{\ell}(X)$ and  $V_{\ell}(Y)$. On the other hand,
there is a natural  embedding of Galois modules \cite[\S
19]{MumfordAV},
$$\Hom(X,Y)\otimes\Q_{\ell}\subset\Hom_{\Q_{\ell}}(V_{\ell}(X),V_{\ell}(Y)),$$
whose image must be a $\Gal(K)$-invariant $\Q_{\ell}$-vector
subspace. It is also clear that
$\Hom_{\Z_{\ell}}(T_{\ell}(X),T_{\ell}(Y))$ is a Galois-invariant
$\Z_{\ell}$-lattice in
$\Hom_{\Q_{\ell}}(V_{\ell}(X),V_{\ell}(Y))$. The equality (3)
gives rise to a natural isomorphism of Galois modules
$$\Hom_{\Z_{\ell}}(T_{\ell}(X),T_{\ell}(Y))\otimes_{\Z_{\ell}}\Z_{\ell}/\ell
\Z_{\ell}=\Hom_{\F_{\ell}}(X_{\ell},Y_{\ell}) \eqno(4).$$

\begin{proof}[Proof of Theorem \ref{hom}]
Let $K(X_{\ell},Y_{\ell})$ be the compositum of the fields
 $K(X_{\ell})$ and $K(Y_{\ell})$. The linear disjointness of
$K(X_{\ell})$ and $K(Y_{\ell})$ means that
$$\Gal(K(X_{\ell},Y_{\ell})/K)= \Gal(K(Y_{\ell})/K)\times
\Gal(K(X_{\ell})/K).$$ Let
$X_{\ell}^*=\Hom_{\F_{\ell}}(X_{\ell},\F_{\ell})$ be the dual of
$X_{\ell}$ and $\bar{\rho}_{n,X,K}^{*}:\Gal(K)
\to\Aut(X_{\ell}^{*})$ the dual of $\bar{\rho}_{n,X,K}$. One may
easily check that
$\ker(\bar{\rho}_{n,X,K}^{*})=\ker(\bar{\rho}_{n,X,K})$ and
therefore we have an isomorphism of the images
$$\tilde{G}^{*}_{\ell,X,K}:=\bar{\rho}_{n,X,K}^{*}(\Gal(K))\cong
\bar{\rho}_{n,X,K}(\Gal(K)))=\tilde{G}_{\ell,X,K}.$$ One may also
easily check that the centralizer of $\Gal(K)$ in
$\End_{\F_{\ell}}(X_{\ell}^{*})$ still coincides with $\F_{\ell}$.
It follows that if $A_1$ is the $\F_{\ell}$-subalgebra in
$\End_{\F_{\ell}}(X_{\ell}^*)$ generated by
$\tilde{G}^{*}_{\ell,X,K}$ then its centralizer in
$\End_{\F_{\ell}}(X_{\ell}^*)$ coincides with $\F_{\ell}$.  Let us
consider the Galois module
$W_1=\Hom_{\F_{\ell}}(X_{\ell},Y_{\ell})=X_{\ell}^{*}\otimes_{\F_{\ell}}
Y_{\ell}$ and denote by $\tau$ the homomorphism
 $\Gal(K)\to\Aut(W_1)$ that defines the Galois module structure on
 $W_1$.
One may easily check that $\tau$ factors through
 $\Gal(K(X_{\ell},Y_{\ell})/K)$ and the image of $\tau$ coincides
 with the image of
$$\tilde{G}^{*}_{\ell,X,K}\times \tilde{G}_{\ell,X,Y}\subset
\Aut(X_{\ell}^{*})\times
\Aut(Y_{\ell})\to\Aut(X_{\ell}^{*}\otimes_{\F_{\ell}}
Y_{\ell})=\Aut(W_1).$$

Let $A_2$ be the $\F_{\ell}$-subalgebra in
$\End_{\F_{\ell}}(Y_{\ell})$ generated by $\tilde{G}_{\ell,Y,K}$.
Recall that the centralizer of $\Gal(K)$ in
$\End_{\F_{\ell}}(Y_{\ell})$ is a field, say $\F$. Clearly,  the
centralizer of $A_2$ in $\End_{\F_{\ell}}(Y_{\ell})$ coincides
with $\F$. One may easily check that the subalgebra of
$\End_{\F_{\ell}}(W_1)$ generated by the image of $\Gal(K)$
coincides with
$$A_1\otimes_{\F_{\ell}}A_2\subset
\End_{\F_{\ell}}(X_{\ell}^*)\otimes_{\F_{\ell}}\End_{\F_{\ell}}(Y_{\ell})=\End_F(X_{\ell}^{*}\otimes_{\F_{\ell}}
Y_{\ell})=\End_{\F_{\ell}}(W_1).$$ It follows from Lemma (10.37)
on p. 252 of \cite{CR} that the centralizer of
$A_1\otimes_{\F_{\ell}}A_2$ in
$\End_F(X_{\ell}^{*}\otimes_{\F_{\ell}} Y_{\ell})$ coincides with
$\F_{\ell}\otimes_{\F_{\ell}}\F=\F$. This implies that the
centralizer of $\Gal(K)$ in
$\End_F(X_{\ell}^{*}\otimes_{\F_{\ell}}
Y_{\ell})=\End_{\F_{\ell}}(W_1)$ is the field $\F$.

Let us consider the $\Q_{\ell}$-vector space
$V_1=\Hom_{\Q_{\ell}}(V_{\ell}(X),V_{\ell}(Y))$ and the free
$\Z_{\ell}$-module $T_1=\Hom_{\Z_{\ell}}(T_{\ell}(X),T_{\ell}(Y))$
provided with the natural structure of Galois modules. Clearly,
$T_1$ is a Galois-stable $\Z_{\ell}$-lattice in $V_1$. By (4),
there is a natural isomorphism of Galois modules $W_1=T_1/\ell
T_1$. Let us denote by $D_1$  the centralizer of $\Gal(K)$ in
$\End_{\Q_{\ell}}(V_1)$. Clearly, $D_1$ is a finite-dimensional
$\Q_{\ell}$-algebra. Therefore in order to prove that $D_1$ is a
division algebra, it suffices to check that $D_1$  has no zero
divisors.

Suppose that $D_1$ has zero divisors, i.e. there are non-zero
$u,v\in D_1$ with $uv=0$. We have $u,v \subset D_1 \subset
\End_{\Q_{\ell}}(V_1)$. Multiplying $u$ and $v$ by proper powers
of $\ell$, we may and will assume that $u (T_1)\subset T_1,
v(T_1)\subset T_1$ but $u(T_1)$ is {\sl not} contained in $\ell
T_1$ and $v(T_1)$ is {\sl not} contained in $\ell T_1$. This means
that $u$ and $v$ induce {\sl non-zero} endomorphisms
$\bar{u},\bar{v}\in \End(W_1)$ that commute with $\Gal(K)$ and
$\bar{u}\bar{v}=0$. Since both $\bar{u}$ and $\bar{v}$ are
non-zero elements of the field $\F$, we get a contradiction that
proves that $D_1$  has no zero divisors and therefore is a
division algebra.

{\sl End of the proof of Theorem} \ref{hom}. We may and will
assume that $K$ is finitely generated over its prime subfield
(replacing $K$ by its suitable subfield). Then the conjecture of
Tate \cite{Tate} (proven by the author in characteristic $>2$
\cite{ZarhinI,ZarhinT}, Faltings in characteristic zero
\cite{Faltings1,Faltings2} and Mori in characteristic $2$
\cite{MB}) asserts that the natural representation of $\Gal(K)$ in
$V_{\ell}(Z)$ is completely reducible for any abelian variety $Z$
over $K$. In particular,   the natural representations of
$\Gal(K)$ in $V_{\ell}(X)$ and $V_{\ell}(Y)$ are completely
reducible. It follows easily that the dual Galois representation
in $\Hom_{\Q_{\ell}}(V_{\ell}(X),\Q_{\ell})$ is also completely
reducible.  Since $\Q_{\ell}$ has characteristic zero, it follows
from a theorem of Chevalley \cite[p. 88]{Chevalley} that the
Galois representation in the tensor product
$\Hom_{\Q_{\ell}}(V_{\ell}(X),\Q_{\ell})\otimes_{\Q_{\ell}}V_{\ell}(Y)=
\Hom_{\Q_{\ell}}(V_{\ell}(X),V_{\ell}(Y))=:V_1$ is completely
reducible. The complete reducibility implies easily that $V_1$ is
an irreducible Galois representation, because the centralizer is a
division algebra. Recall that $\Hom(X,Y)\otimes\Q_{\ell}$ is a
Galois-invariant subspace in
$\Hom_{\Q_{\ell}}(V_{\ell}(X),V_{\ell}(Y))=V_1$. The
irreducibility of $V_1$ implies that either
$\Hom(X,Y)\otimes\Q_{\ell}=0$ or $\Hom(X,Y)\otimes\Q_{\ell}=V_1$.

 If $\Hom(X,Y)\otimes\Q_{\ell}=0$ then $\Hom(X,Y)=0$ and therefore
$\Hom(Y,X)=0$.

If $\Hom(X,Y)\otimes\Q_{\ell}=V_1$ then the rank of the free
commutative group $\Hom(X,Y)$ coincides with the dimension of the
$\Q_{\ell}$-vector space $V_1$. Clearly, $V_1$ has dimension
$4\dim(X)\dim(Y)$. It is proven in Proposition 3.3 of
\cite{ZarhinSh} that if $A$ and $B$ are abelian varieties over an
algebraically closed field ${\mathcal{K}}$ and the rank of
$\Hom(A,B)$ equals $4\dim(A)\dim(B)$ then
$\fchar({\mathcal{K}})>0$ and both $A$ and $B$ are supersingular
abelian varieties. Applying this result to $X$ and $Y$, we
conclude that $\fchar(K)=\fchar(K_a)>0$ and both $X$ and $Y$ are
supersingular abelian varieties.
\end{proof}

\section{Hyperelliptic jacobians}
\label{jac} In this section we deal with the case of $\ell=2$.
 Suppose that $\fchar(K)\ne 2$. Let $f(x)\in K[x]$ be a polynomial
of degree $n\ge 3$ without multiple roots. Let $\RR_f\subset K_a $
be the set of roots of $f$. Clearly, $\RR_f$ consists of $n$
elements. Let $K(\RR_f)\subset K_a$ be the splitting field of $f$.
Clearly, $K(\RR_f)/K$ is a Galois extension and we write $\Gal(f)$
for its Galois group $\Gal(K(\RR_f)/K)$. By definition,
$\Gal(K(\RR_f)/K)$ permutes elements of  $\RR_f$; further we
identify $\Gal(f)$ with the corresponding subgroup of
$\Perm(\RR_f)$ where $\Perm(\RR_f)$ is the group of permutations
of $\RR_f$.

We write $\F_2^{\RR_f}$ for the $n$-dimensional $\F_2$-vector
space of maps $h:\RR_f \to \F_2$. The space $\F_2^{\RR_f}$ is
provided with a natural action of $\Perm(\RR_f)$ defined as
follows. Each $s \in \Perm(\RR_f)$ sends a map
 $h:\RR_f\to \F_2$ to  $sh:\alpha \mapsto h(s^{-1}(\alpha))$. The permutation module $\F_2^{\RR_f}$
contains the $\Perm(\RR_f)$-stable hyperplane
$$(\F_2^{\RR_f})^0=
\{h:\RR_f\to \F_2\mid\sum_{\alpha\in \RR_f}h(\alpha)=0\}$$ and the
$\Perm(\RR_f)$-invariant line $\F_2 \cdot 1_{\RR_f}$ where
$1_{\RR_f}$ is the constant function $1$. Clearly,
$(\F_2^{\RR_f})^0$ contains $\F_2 \cdot 1_{\RR_f}$ if and only if
$n$ is even.

 If  $n$ is even then
let us define the $\Gal(f)$-module
$Q_{\RR_f}:=(\F_2^{\RR_f})^0/(\F_2 \cdot 1_{\RR_f})$. If $n$ is
odd then let us put $Q_{\RR_f}:=(\F_2^{\RR_f})^0$. If $n \ne 4$
the natural representation of $\Gal(f)$ is faithful, because in
this case the natural homomorphism
$\Perm(\RR_f)\to\Aut_{\F_2}(Q_{\RR_f})$ is injective.

\begin{rem}
\label{klemm}
 It is known \cite[Satz 4]{Klemm},  that
$\End_{\Gal(f)}(Q_{\RR_f})=\F_2$ if either $n$ is odd and
$\Gal(f)$ acts doubly transitively on $\RR_f$ or $n$ is even and
$\Gal(f)$ acts  $3$-transitively on $\RR_f$.
\end{rem}

The canonical surjection $\Gal(K)\twoheadrightarrow
\Gal(K(\RR_f)/K)=\Gal(f)$ provides $Q_{\RR_f}$ with a natural
structure of $\Gal(K)$-module. Let $C_f$ be the hyperelliptic
curve $y^2=f(x)$ and $J(C_F)$ its jacobian. It is well-known that
$J(C_F)$ is a $\left[\frac{n-1}{2}\right]$-dimensional abelian
variety defined over $K$. It is also well-known that the
$\Gal(K)$-modules $J(C_f)_2$ and $Q_{\RR_f}$ are isomorphic (see
for instance \cite{Poonen,SPoonen,ZarhinTexel}). It follows that
if $n\ne 4$ then
$$\Gal(f)=\tilde{G}_{2,J(C_f)}.$$
It follows from Remark \ref{klemm} that if either $n$ is odd and
$\Gal(f)$ acts doubly transitively on $\RR_f$ or $n$ is even and
$\Gal(f)$ acts  $3$-transitively on $\RR_f$ then
$$\End_{\tilde{G}_{2,J(C_f)}}(J(C_f)_2))=\F_2.$$
It is also clear that $K(J(C_f)_2))\subset K(\RR_f)$. (The
equality holds if  $n\ne 4$.)

The next assertion follows immediately from Theorem \ref{center2},
Corollaries \ref{odd} and \ref{int} (applied to
$X=J(C_f),\ell=2,\GG=\Gal(f)$).

\begin{thm}
\label{double}
 Let $K$ be a field of characteristic different from
$2$, let $n\ge 5$ be an integer, $g=\left[\frac{n-1}{2}\right]$
and $f(x)\in K[x]$ a polynomial of degree $n$. Suppose that either
$n$ is odd and $\Gal(f)$ acts doubly transitively on $\RR_f$ or
$n$ is even and $\Gal(f)$ acts  $3$-transitively on $\RR_f$.
Assume also that $\Gal(f)$ is a simple nonabelian group that does
not contain a subgroup of index dividing $g$ except $\Gal(f)$
itself. If $g$ is odd then $\End^0(J(C_f))$ enjoys one of the
following properties:
\begin{itemize}
\item[(i)]
$\End^0(J(C_f))$ is isomorphic to the matrix algebra $\M_d(\Q)$
where $d$ divides $g$. If $d>1$ there exist a finite perfect group
$\Pi\subset\GL(d,\Z)$ and a surjective homomorphism
$\Pi\twoheadrightarrow\Gal(f)$ such that every prime dividing
$\#(\Pi)$ also divides $\#(\Gal(f))$.
\item[(ii)]
$p:=\fchar(K)$ is a prime dividing $\#(\Gal(f))$ and
$\End^0(J(C_f))$ is isomorphic to the matrix algebra $\M_d(\H_p)$
where $d>1$ divides $g$.
\end{itemize}
\end{thm}

\begin{ex}
\label{a5} Suppose that $n=5$ and $\Gal(f)$ is the alternating
group $\A_5$ acting doubly transitively on $\RR_f$. Clearly, $g=2$
and $\Gal(f)$ is a simple nonabelian group without subgroups of
index $2$.  Applying Theorem \ref{double}, we conclude that
$\End^0(J(C_f))$ is either $\Q$ or $\M_2(\Q)$ or $\M_2(\H)$ where
$\H$ is a quaternion $\Q$-algebra unramified outside
$\{\infty,2,3,5\}$; in addition $\H\cong \H_p$ if
$p:=\fchar(K)>0$. Suppose that $\End(J(C_f))\ne \Z$ and therefore
$\End^0(J(C_f))\ne \Q$. If $\End^0(J(C_f))=\M_2(\Q)$ then
$\GL(2,\Q)=\M_2(\Q)^{*}$ contains a finite group, whose order
divides $5$, which is not the case. This implies that
$\End^0(J(C_f))=\M_2(\H)$. This means that $J(C_f)$ is
supersingular and therefore  $p:=\fchar(K)>0$. This implies that
$p=3$ or $p=5$.

We conclude that either $\End(J(C_f))= \Z$ or $\fchar(K)\in
\{3,5\}$ and $J(C_f)$ is a supersingular abelian varietiy. In
fact, it is known \cite{ZarhinBSMF} that if $\fchar(K)=5$ then
$\End(J(C_f))= \Z$. On the other hand, one may find a
supersingular $J(C_f)$ in characteristic $3$ \cite{ZarhinBSMF}.
\end{ex}

Example \ref{a5} is a special case of the following general result
proven by the author \cite{ZarhinMRL,ZarhinMMJ,ZarhinBSMF}. {\sl
Suppose that $n\ge 5$ and $\Gal(f)$ is the alternating group
$\A_n$ acting on $\RR_f$. If $\fchar(K)=3$ we assume additionally
that $n\ge 7$. Then $\End(J(C_f))= \Z$}.

We refer the reader to
~\cite{Mori1,Mori2,Katz1,Katz2,Masser,KS,ZarhinMRL,ZarhinMRL2,ZarhinMMJ,ZarhinP,ZarhinPAMS,ZarhinSt}
for a discussion of other known results about, and examples of,
hyperelliptic jacobians without complex multiplication.

\begin{cor}
\label{l27} Suppose that $n=7$ and
$\Gal(f)=\SL_3(\F_2)\cong\PSL_2(\F_7)$ acts doubly transitively on
$\RR_f$. Then $\End^0(J(C_f))=\Q$ and therefore $\End(J(C_f))=\Z$.
\end{cor}
\begin{proof}
We have $g=\dim(J(C_f))=3$. Since $\PSL_2(\F_7)$ is a simple
nonabelian group it does not contain a subgroup of index $3$. So,
we may apply Theorem \ref{double}. We obtain that if
$\End^0(J(C_f))\ne\Q$ then either $\End^0(J(C_f))=\M_3(\Q)$ and
there exist a finite perfect group $\Pi\subset\GL(3,\Z)$ and a
surjective homomorphism
$\Pi\twoheadrightarrow\Gal(f)=\PSL_2(\F_7)$ or
$\End^0(J(C_f))=\M_3(\H_p)$ where $p=\fchar(K)$ is either $3$ or
$7$. The  case of $\End^0(J(C_f))=\M_3(\H_p)$ means that $J(C_f)$
is supersingular,  which is not true \cite[Th. 3.1]{ZarhinBSMF}.
Hence $\End^0(J(C_f))=\M_3(\Q)$ and $\GL(3,\Z)$ contains a finite
group, whose order is divisible by $7$. It follows that
$\GL(3,\Z)$ contains an element of order  $7$, which is not true.
The obtained contradiction proves that $\End^0(J(C_f))=\Q$ and
therefore $\End(J(C_f))=\Z$.
\end{proof}

\begin{cor}
\label{l211} Suppose that $n=11$ and $\Gal(f)=\PSL_2(\F_{11})$
acts doubly transitively on $\RR_f$. Then $\End^0(J(C_f))=\Q$ and
therefore $\End(J(C_f))=\Z$.
\end{cor}
\begin{proof}
We have $g=\dim(J(C_f))=5$. It is known \cite{Atlas} that
$\PSL_2(\F_{11})$ is a simple nonabelian subgroup  not containing
a subgroup of index $5$. So, we may apply Theorem \ref{double}. We
obtain that if $\End^0(J(C_f))\ne\Q$ then either
$\End^0(J(C_f))=\M_5(\Q)$ and there exist a finite perfect group
$\Pi\subset\GL(5,\Z)$ and a surjective homomorphism
$\Pi\twoheadrightarrow\Gal(f)=\PSL_2(\F_{11})$ or
$\End^0(J(C_f))=\M_5(\H_p)$ where $p=\fchar(K)$ is either $3$ or
$5$ or $11$.

Assume that  $\End^0(J(C_f))=\M_5(\Q)$. Then $\GL(5,\Z)$ contains
a finite group,  whose order is divisible by $11$. It follows that
$\GL(5,\Z)$ contains an element of order $11$, which is not true.
Hence $\End^0(J(C_f))\ne\M_5(\Q)$.

Assume that $\End^0(J(C_f))=\M_5(\H_p)$ where $p$ is either $3$ or
$5$ or $11$. This implies that  $J(C_f)$ is a supersingular
abelian variety.

Notice that every homomorphism from simple $\PSL_{2}(\F_{11})$ to
$\GL(4,\F_2)$ is trivial, because $11$ divides
$\#(\PSL_{2}(\F_{11}))$  but $\#(\GL(4,\F_2))$ is {\sl not}
divisible by $11$. Since $4=g-1$, it follows from Theorem 3.3 of
\cite{ZarhinBSMF} (applied to $g=5,
X=J(C_f),G=\Gal(f)=\PSL_{2}(\F_{11})$) that there exists a central
extension $\pi_1:G_1 \to \PSL_{2}(\F_{11})$ such that $G_1$ is
perfect,  $\ker(\pi_1)$ is a cyclic group of order $1$ or $2$ and
$\M_5(\H_p)$ is a direct summand of the group $\Q$-algebra
$\Q[G_1]$. It follows easily that  $G_1=\PSL_{2}(\F_{11})$ or
$\SL_{2}(\F_{11})$. It is known \cite{J,F} that
$\Q[\PSL_{2}(\F_{11})]$ is a direct sum of matrix algebras over
fields. Hence $G_1=\SL_{2}(\F_{11})$ and the direct summand
$\M_5(\H_p)$ corresponds to a faithful ordinary irreducible
character $\chi$ of $\SL_{2}(\F_{11})$ with degree $10$ and
$\Q(\chi)=\Q$. This implies that in notations of \cite[\S 38]{DA},
 $\chi=\theta_j$ where $j$ is an odd integer such that $1\le j
\le \frac{11-1}{2}=5$ and either $6j$ is divisible by $11+1=12$ or
$4j$ is divisible by $12$ (\cite{F}, Th. 6.2 on p. 285). This
implies that $j=3$ and $\chi=\theta_3$. However, the direct
summand attached to $\theta_3$ is ramified at $2$ (\cite[the case
(c) on p. 4]{J};  \cite[theorem  6.1(iii) on p. 284]{F}). Since $p
\ne 2$, we get a contradiction which
  proves that $J(C_f)$ is not
supersingular. This implies that $\End^0(J(C_f))=\Q$ and therefore
$\End(J(C_f))=\Z$.
\end{proof}

\begin{cor}
\label{m12} Suppose that $n=12$ and $\Gal(f)$ is the Mathieu group
$\M_{12}$ acting  $3$-transitively on $\RR_f$. Then
$\End(J(C_f))=\Z$.
\end{cor}
\begin{proof}
Let $\alpha$ be a root of $f(x)$ and $K_1=K(\alpha)$. Clearly, the
stabilizer of $\alpha$ in $\Gal(f)=\M_{12}$ is $\PSL_2(\F_{11})$
acting doubly transitively on the roots of
$f_1(x)=\frac{f(x)}{x-\alpha}\in K_1[x]$. Let us put
$h(x)=f_1(x+\alpha)\in K_1[x],  h(x)=x^{11} h(1/x)\in K_1[x]$.
Clearly, $\deg(h_1)=11$ and  $\Gal(h_1)=\PSL_2(\F_{11})$ acts
doubly transitively on the roots of $h_1$. By Corollary
\ref{l211}, $\End(J(C_{h_1}))=\Z$. On the other hand, the standard
substitution $x_1=1/(x-\alpha),  y_1=y/(x-\alpha)^6$ establishes a
birational isomorphism between $C_f$ and $C_{h_1}:y_1^2=h_1(x_1)$.
This implies that $J(C_f)\cong J(C_{h_1})$ and therefore
$\End(J(C_f))=\Z$.
\end{proof}

In characteristic zero the assertions of Corollaries \ref{l27},
\ref{l211} and \ref{m12} were earlier proven in
\cite{ZarhinBSMF,ZarhinTexel}.

\begin{cor}
\label{mat} Suppose that $\deg(f)=n$ where $n=22,23$ or $24$ and
$\Gal(f)$ is the corresponding (at least) $3$-transitive Mathieu
group $\MM_n\subset\Perm(\RR_f)\cong \Sn$. Then $\End(J(C_f))=\Z$.
\end{cor}
\begin{proof}
First, assume that $n=23$ or $24$. We have $g=\dim(J(C_f))=11$. It
is known that both $\MM_{23}$ and $\MM_{24}$ do not contain a
subgroup of index $11$ \cite{Atlas}. So, we may apply Theorem
\ref{double} and obtain that if $\End(J(C_f)\ne\Z$ then
$\End^0(J(C_f))\ne\Q$ and one of the following conditions holds:

\begin{itemize}
\item[(i)]
$\End^0(J(C_f))=\M_{11}(\Q)$ and there exist a finite perfect
group $\Pi\subset\GL(11,\Z)$ and a surjective homomorphism
$\Pi\twoheadrightarrow\Gal(f)=\MM_n$;
\item[(ii)]
$p=\fchar(K)\in\{3,5,7,11,23\}$ and
$\End^0(J(C_f))=\M_{11}(\H_p)$.
\end{itemize}

Assume that the condition (i) holds. Then
$\End^0(J(C_f))=\M_{11}(\Q)$ and $\GL(11,\Z)$ contains a finite
group, whose  order is divisible by $23$. It follows that
$\GL(11,\Z)$ contains an element of order  $23$, which is not
true. The obtained contradiction proves that the condition (i) is
not fulfilled.

Hence  the condition (ii) holds. Then
$p=\fchar(K)\in\{3,5,7,11,23\}$ and  there exist a finite perfect
subgroup $\Pi\subset \End^0(J(C_f))^{*}=\GL(11,\H_p)$ and a
surjective homomorphism $\pi:\Pi\twoheadrightarrow\MM_n$.
Replacing $\Pi$ by a suitable subgroup, we may and will assume
that no proper subgroup of $\Pi$ maps onto $\MM_{n}$. By tensoring
$\H_p$ to the field of complex numbers (over $\Q$), we obtain an
embedding
$$\Pi \subset \GL(11,\H_p)\subset\GL(22,\C).$$
In particular, the (perfect) group $\Pi$ admits a non-trivial
projective  $22$-dimensional representation over $\C$. Recall that
$\MM_{n}$ has Schur's multiplier $1$ (since $n=23$ or $24$)
\cite{Atlas} and therefore all its projective representations are
(obtained from) linear representations. Also, all nontrivial
linear representations of $\MM_{24}$ have dimension $\ge 23$,
because the smallest dimension of a nontrivial linear
representation of $\MM_{24}$ is $23$. It follows from results of
Feit--Tits \cite{FT} that $\Pi$ cannot have  a non-trivial
projective representation of dimension $<23$. This implies that $n
\ne 24$, i.e. $n=23$.

Recall that $22$ is the smallest possible dimension of a
nontrivial representation of $\MM_{23}$ in characteristic zero,
because its every irreducible representation in characteristic
zero has dimension $\ge 22$ \cite{Atlas}. It follows from a
theorem of Feit--Tits (\cite{FT}, pp. 1 and \S 4; see also
\cite{KL}) that the projective representation
$$\Pi \to \GL(11,\H_p)/\Q^{*}\subset\GL(22,\C)/\C^{*}$$
factors through $\ker(\pi)$. This means that $\ker(\pi)$ lies in
$\Q^{*}$ and therefore $\Pi$ is a central extension of $\MM_{23}$.
Now the perfectness of $\Pi$ implies that $\pi$ is an isomorphism,
i.e. $\Pi\cong\MM_{23}$.

Let us consider the natural homomorphism $\Q[\MM_{23}]\cong
\Q[\Pi] \to \M_{11}(\H_p)$ induced by the inclusion $\Pi\subset
\M_{11}(\H_p)^{*}$. It is surjective, because otherwise one may
construct a (complex) nontrivial representation of $\MM_{23}$ of
dimension $<22$. This implies that $\M_{11}(\H_p)$ is isomorphic
to a direct summand of $\Q[\MM_{23}]$. But this is not true, since
Schur indices of all irreducible representations of $\MM_{23}$ are
equal to $1$ \cite[\S 7]{F} and therefore $\Q[\MM_{23}]$ splits
into a direct sum of matrix algebras over fields. The obtained
contradiction proves that the condition (ii) is not fulfilled. So,
$\End(J(C_f))=\Z$.

Now let $n=22$.  Then $g=10$. It is known that $\MM_{22}$ is a
simple nonabelian group
 not containing a subgroup of index $10$ \cite{Atlas}. Let us assume that
 $\End^0(J(C_f))\ne\Q$.
 Applying Theorem
\ref{center2}, we conclude that there exists a positive integer
$d$ dividing $10$ such that either $d>1$ and
$\End^0(J(C_f))=\M_d(\Q)$ or $\End^0(J(C_f))=\M_d(\H)$ where $\H$
is a quaternion $\Q$-algebra unramified outside $\infty$ and the
prime divisors of $\#(\MM_{22})$. In addition, there exist a
finite perfect subgroup $\Pi\subset \End^0(J(C_f))^{*}$ and a
surjective homomorphism $\pi:\Pi\twoheadrightarrow\MM_{22}$.
Replacing $\Pi$ by a suitable subgroup, we may and will assume
(without losing the perfectness) that no proper subgroup of $\Pi$
maps onto $\MM_{n}$.

 By Lemma 3.13 on pp. 200--201 of \cite{ZarhinP}, every
homomorphism from $\Pi$ to $\PSL(10,\R)$ is trivial. The
perfectness of $\Pi$ implies that every homomorphism from $\Pi$ to
$\PGL(10,\R)$ is trivial. Since $\M_d(\Q)^*=\GL(d,\Q)\subset
\GL(10,\R)$, we conclude that $\End^0(J(C_f))\ne \M_d(\Q)$ and
therefore $\End^0(J(C_f))=\M_d(\H)$.

If $d=10$ then $p:=\fchar(K)>0$ and  $J(C_f)$ is a supersingular
abelian variety.

 Assume that $d \ne 10$, i.e. $d=1,2$ or $5$. If $H$ is unramified at $\infty$
 then there exists an embedding $\H \hookrightarrow \M_2(\R)$.
 This gives us the embeddings
$$\Pi\subset \M_d(\H)^{*}\hookrightarrow \M_{2d}(\R)^{*}=\GL(2d,\R)\subset\GL(10,\R)$$
and therefore there is a nontrivial homomorphism from $\Pi$ to
$\PGL(10,\R)$. The obtained contradiction proves that $\H$ is
ramified at $\infty$.

There exists an embedding $\H\hookrightarrow \M_4(\Q)\subset
\M_4(\R)$. This implies that if $d=1$ or $2$ then there are
embeddings
$$\Pi\subset \M_d(\H)^{*}\hookrightarrow \M_{4d}(\R)^{*}=\GL(4d,\R)\subset\GL(10,\R)$$
and therefore there is a nontrivial homomorphism from $\Pi$ to
$\PGL(10,\R)$. The obtained contradiction proves that $d=5$. This
means that there exists an abelian surface $Y$ over $K_a$ such
that $J(C_f)$ is isogenous to $Y^5$ and $\End^0(Y)=\H$. However,
there do not exist abelian surfaces, whose endomorphism algebra is
a definite quaternion algebra over $\Q$. This result is well-known
in characteristic zero (see, for instance \cite{OZ}); the positive
characteristic  case was done by Oort \cite[Lemma 4.5 on p.
490]{Oort}. Hence $d\ne 5$. This implies that $d=10$ and $J(C_f)$
is a supersingular abelian variety.

 Since
$\MM_{22}$ is a simple group and $11\mid\#(\MM_{22})$, every
homomorphism from $\MM_{22}$ to $\GL(9,\F_2)$ is trivial, because
$\#(\GL(9,\F_2))$ is not divisible by $11$. Since $9=g-1$, it
follows from Theorem 3.3 of \cite{ZarhinBSMF} (applied to $g=10,
X=J(C_f),G=\Gal(f)=\MM_{22}$) that there exists a central
extension $\pi_1:G_1 \to \MM_{22}$ such that $G_1$ is perfect,
$\ker(\pi_1)$ is a cyclic group of order $1$ or $2$ and there
exists a faithful $20$-dimensional absolutely irreducible
representation of $G_1$ in characteristic zero. However, such a
central extension with $20$-dimensional irreducible representation
does not exist \cite{Atlas}.
\end{proof}

Combining Corollary \ref{mat} with previous author's results
\cite{ZarhinTexel,ZarhinMMJ} concerning small Mathieu groups, we
obtain the following statement.

\begin{thm}
\label{m} Suppose that $n\in\{11,12,22,23,24\}$ and $\Gal(f)$ is
the corresponding  Mathieu group $\MM_n\subset\Perm(\RR_f)\cong
\Sn$. Then $\End(J(C_f))=\Z$.
\end{thm}

In characteristic zero  the assertion of Theorem \ref{m} was
earlier proven in \cite{ZarhinTexel,ZarhinP}.

\begin{thm}
Suppose that $n=15$ and $\Gal(f)$ is the alternating group $\A_7$
acting doubly transitively on $\RR_f$. Then either
$\End(J(C_f))=\Z$ or $J(C_f)$ is isogenous over $K_a$ to a product
of elliptic curves.
\end{thm}

\begin{proof}
We have $g=7$. Unfortunately, $\A_7$ has a subgroup of index $7$.
However, $\A_7$ is simple nonabelian and does not have a normal
subgroup of index $7$. Applying Theorem \ref{center2} to
$X=J(C_f),g=7, \ell=2,\GG=\Gal(f)=\A_7$,
 we obtain that either $J(C_f)$ is isogenous to a product of
 elliptic curves (case (a))
  or $\End^0(J(C_f))$ is a central simple
 $\Q$-algebra (case (b)). If $\End^0(J(C_f))$ is a matrix algebra over $\Q$
 then either $\End^0(J(C_f))=\Q$ (i.e., $\End(J(C_f))=\Z$) or $\End^0(J(C_f))=\M_7(\Q)$
 (i.e., $J(C_f)$ is isogenous to the $7$th power of an elliptic curve
 without complex multiplication).

 If the central simple
 $\Q$-algebra $\End^0(J(C_f))$ is not a matrix algebra over $\Q$ then there
 exists a quaternion $\Q$-algebra $\H$ such that either
 $\End^0(J(C_f))=\H$ or $\End^0(J(C_f))=\M_7(\H)$. If
 $\End^0(J(C_f))=\M_7(\H)$ then $J(C_f)$ is a supersingular abelian
 variety and therefore is isogenous to a product of elliptic
 curves.

 Let us assume that $\End^0(J(C_f))=\H$. We need to arrive to a
 contradiction. Since $7=\dim(J(C_f))$ is odd,
 $p=\fchar(K)>0$. The same arguments as in the proof of Corollary
 \ref{odd} tell us that $\H=\H_p$. By Theorem \ref{center2}(b3), there
 exist a perfect finite group $\Pi\subset\End^0(J(C_f))^*=\H_p^{*}$ and a
 surjective homomorphism $\Pi\twoheadrightarrow\A_7$. But Lemma \ref{Hp} asserts
 that every finite subgroup in $\H_p^{*}$ is solvable. The obtained
 contradiction proves that $\End^0(J(C_f))\ne\H_p$.
\end{proof}

\begin{thm}
\label{l2q} Suppose that $n=q+1$ where $q\ge 5$ is a prime  power
that is congruent to $\pm 3$ modulo $8$. Suppose that
$\Gal(f)=\PSL_2(\F_q)$ acts doubly transitively on $\RR_f$ (where
$\RR_f$ is identified with the projective line $\P^1(\F_q)$). Then
$\End^0(J(C_f))$ is a simple $\Q$-algebra, i.e. $J(C_f)$ is either
absolutely simple or isogenous to a power of an absolutely simple
abelian variety.
\end{thm}

\begin{proof}
Since $n=q+1$ is even, $g=\frac{q-1}{2}$. It is known
\cite{Mortimer} that the $\Gal(f)=\PSL_2(\F_q)$-module $Q_{\RR_f}$
is simple and the centralizer of $\PSL_2(\F_q)$ in
$\End_{\F_2}(Q_{\RR_f})$ is the field $\F_4$. On the other hand,
$\PSL_2(\F_q)$ is a simple nonabelian group: we need to inspect
its subgroups. The following statement will be proven later in
this section.

\begin{lem}
\label{sub} Let $q\ge 5$ be a  power of an odd prime. Then
$\PSL_2(\F_q)$ does not contain a subgroup of index dividing
$\frac{q-1}{2}$ except $\PSL_2(\F_q)$ itself.
\end{lem}

 Recall that
$\tilde{G}_{2,J(C_f)}=\Gal(f)=\PSL_2(\F_q)$. Now Theorem \ref{l2q}
follows readily from Theorem \ref{center1} combined with Lemma
\ref{sub}.
\end{proof}

\begin{proof}[Proof of Lemma \ref{sub}]
Since $\PSL_2(\F_q)$ is a simple nonabelian subgroup, it does not
contain a subgroup of index $\le 4$ except $\PSL_2(\F_q)$ itself.
This implies that in the course of the proof we may assume that
$\frac{q-1}{2}\ge 5$, i.e., $q\ge 11$.

 Recall that $\#(\PSL_2(\F_q))=(q+1)q(q-1)/2$. Let $H \ne
\PSL_2(\F_q)$ be a subgroup in $\PSL_2(\F_q)$. The list of
subgroups in $\PSL_2(\F_q)$ given in \cite[theorem 6.25 on p.
412]{Suzuki} tells us that $\#(H)$ divides either $q\pm 1$ or
$q(q-1)/2$ or $60$ or $(b+1)b(b-1)$ where $b<q$ is a positive
integer such that $q$ is an integral power of $b$. This implies
that if the index of $H$ is a divisor of $\frac{q-1}{2}$ then
either

\begin{enumerate}
\item
$(q+1)q$ divides  $60$

or
\item
$\frac{(q+1)q(q-1)}{2} \le \frac{q-1}{2}
(\sqrt{q}+1)\sqrt{q}(\sqrt{q}-1)=\frac{q-1}{2}(q-1)\sqrt{q}$.
\end{enumerate}

In the case (1) we have $q=5$  which contradicts our assumption
that $q\ge 11$. So, the case (2) holds. Clearly, $(q+1)\sqrt{q}
\le (q-1)$ which is obviously not true.
\end{proof}

\begin{thm}
\label{homjac} Let $K$ be a field of characteristic different from
$2$. Suppose that $f(x)$ and $h(x)$ are polynomials in $K[x]$
enjoying the following properties:
\begin{itemize}
\item[(i)]
$\deg(f)\ge 3$ and the Galois group $\Gal(f)$ acts doubly
transitively on the set $\RR_f$ of roots of $f$. If $\deg(f)$ is
even then this action is $3$-transitive;
\item[(ii)]
 $\deg(h)\ge 3$ and the Galois group $\Gal(h)$ acts doubly
transitively on the set $\RR_h$ of roots of $h$. If $\deg(h)$ is
even then this action is $3$-transitive;
\item[(iii)]
The splitting fields $K(\RR_f)$ of $f$ and $K(\RR_h)$ of $h$ are
linearly disjoint over $K$.
\end{itemize}
Let $J(C_f)$ be the jacobian of the hyperelliptic curve
$C_f:y^2=f(x)$ and $J(C_h)$ be the jacobian of the hyperelliptic
curve $C_h:y^2=h(x)$. Then either $\Hom(J(C_f),J(C_h))=0,
\Hom(J(C_h),J(C_f))=0$ or $\fchar(K)>0$ and both $J(C_f)$ and
$J(C_h)$ are supersingular abelian varieties.
\end{thm}

\begin{proof}
Let us put $X=J(C_f), Y=J(C_h)$. The transitivity properties imply
that $\End_{\tilde{G}_{2,X}}(X_2)=\F_2$ and
$\End_{\tilde{G}_{2,Y}}(Y_2)=\F_2$. The linear disjointness of
$K(\RR_f)$ and $K(\RR_h)$ implies that the fields
$K(X_2)=K((J(C_f)_2)\subset K(\RR_f)$ and
$K(Y_2)=K((J(C_h)_2)\subset K(\RR_h)$ are also linearly disjoint
over $K$. Now the assertion follows readily from Theorem \ref{hom}
with $\ell=2$.
\end{proof}

\section{Abelian varieties with multiplications}
\label{mult}

Let $E$ be a number field. Let $(X, i)$ be a pair consisting of an
abelian variety $X$ of positive dimension over $K_a$ and an
embedding $i:E \hookrightarrow  \End^0(X)$. Here $1\in E$ must go
to $1_X$. It is well known \cite{Ribet2} that the degree $[E:\Q]$
divides $2\dim(X)$, i.e.
$$d=d_X:=\frac{2\dim(X)}{[E:\Q]}$$
is a positive integer. Let us denote by $\End^0(X,i)$ the
centralizer of $i(E)$ in $\End^0(X)$. The image $i(E)$ lies in the
center of the finite-dimensional $\Q$-algebra $\End^0(X,i)$. It
follows that $\End^0(X,i)$ carries a natural structure of
finite-dimensional $E$-algebra. If $Y$ is (possibly) another
abelian variety over $K_a$ and $j:E \hookrightarrow \End^0(Y)$ is
an embedding that sends $1$ to  $1_Y$ then we write
$$\Hom^0((X,i),(Y,j))=\{u \in \Hom^0(X,Y)\mid ui(c)=j(c)u \quad
\forall c\in E\}.$$ Clearly, $\End^0(X,i)=\Hom^0((X,i),(X,i))$. If
$m$ is a positive integer then we write $i^{(m)}$ for the
composition $E\hookrightarrow \End^0(X)\subset \End^0(X^m)$ of $i$
and the diagonal inclusion $\End^0(X)\subset
\End^0(X^m)=\M_m(\End^0(X))$. We have
$$\End^0(X^m,i^{(m)})=\M_m(\End^0(X,i))\subset\M_m(\End^0(X))=\End^0(X^m).$$

\begin{rem}
\label{ss} The $E$-algebra $\End^0(X,i)$ is semisimple. Indeed, in
notations of Remark \ref{split} $\End^0(X)= \prod_{s\in\I} D_s$
where all  $D_s=\End^0(X_s)$ are simple $\Q$-algebras. If
$\pr_s:\End^0(X) \twoheadrightarrow D_s$ is the corresponding
projection map and
  $D_{s,E}$ is the centralizer of $\pr_s i(E)$ in $D_s$ then one may easily check that
$\End^0(X,i)=\prod_{s\in\I} D_{s,E}$. Clearly, $\pr_s i(E)\cong E$
is a simple $\Q$-algebra. It follows from Theorem 4.3.2 on p. 104
of \cite{Herstein} that $D_{s,E}$ is also a {\sl simple}
$\Q$-algebra. This implies  that $D_{s,E}$ is a {\sl simple}
$E$-algebra and therefore $\End^0(X,i)$ is a semisimple
$E$-algebra. We write $i_s$ for the composition $ \pr_s i: E
\hookrightarrow \End^0(X)\twoheadrightarrow D_{s} \cong
\End^0(X_s)$. Clearly, $D_{s,E}=\End^0(X_s,i_s)$ and
$$\End^0(X,i)=\prod_{s\in\I} \End^0(X_s,i_s) \eqno(5) .$$
It follows that $\End^0(X,i)$ is a simple $E$-algebra if and only
if $\End^0(X)$ is a simple $\Q$-algebra, i.e., $X$ is isogenous to
a self-product of (absolutely) simple abelian variety.
\end{rem}

\begin{thm}
\label{maxE}
\begin{itemize}
\item[(i)]
$\dim_E(\End^0((X,i)) \le \frac{4\cdot\dim(X)^2}{[E:\Q]^2}$;
\item[(ii)]
Suppose that
 $\dim_E(\End^0((X,i)) =
\frac{4\cdot\dim(X)^2}{[E:\Q]^2}$.
 Then:
\begin{itemize}
\item[(a)]
 $X$ is isogenous to a self-product of an
(absolutely) simple abelian variety. Also $\End^0((X,i)$ is a
central simple $E$-algebra, i.e., $E$ coincides with the center of
$\End^0((X,i)$. In addition, $X$ is an abelian variety of CM-type.
\item[(b)]
There exist an abelian variety $Z$, a positive integer $m$, an
isogeny $\psi: Z^m \to X$ and an embedding $k: E \hookrightarrow
\End^0(Z)$ that sends $1$ to $1_Z$ such that:

\begin{enumerate}
\item
$\End^0(Z,k)$ is a central division algebra over $E$ of dimension
$\left(\frac{2\dim(Z)}{[E:\Q]}\right)^2$ and $\psi \in
\Hom^0((Z^r,k^{(m)}),(X,i))$.
\item
If $\fchar(K_a)=0$ then $E$ contains a CM subfield and
$2\dim(Z)=[E:\Q]$. In particular, $[E:\Q]$ is even.
\item
If $E$ does not contain a CM-field (e.g., $E$ is a totally real
number field) then $\fchar(K_a)>0$ and $X$ is a supersingular
abelian variety.
\end{enumerate}

\end{itemize}
\end{itemize}
\end{thm}

\begin{proof}
Recall that $d=2\dim(X)/[E:\Q]$. First, assume that $X$ is
isogenous to a self-product of an absolutely simple abelian
variety, i.e., $\End^0(X,i)$ is a simple $E$-algebra. We need to
prove that
$$N:=\dim_E(\End^0(X,i))\le d^2.$$
Let $C$ be the center of $\End^0(X)$.
 Let $E'$ be the center of
$\End^0(X,i)$. Clearly,
$$C \subset E'\subset \End^0(X,i)\subset\End^0(X).$$

Let us put $e=[E':E]$. Then $\End^0(X,i)$ is a {\sl central}
simple $E'$-algebra of dimension $N/e$. Then there exists a
central division $E'$-algebra $D$ such that $\End^0(X,i)$ is
isomorphic to the matrix algebra $\M_m(D)$ of size $m$ for some
positive integer $m$. Dimension arguments imply that
$$m^2\dim_{E'}(D)=\frac{N}{e},\quad \dim_{E'}(D)=\frac{N}{e m^2}.$$
Since $\dim_{E'}(D)$ is a square,
$$\frac{N}{e}=N_1^2, \quad N=e N_1^2, \quad \dim_{E'}(D)=\left(\frac{N_1}{m}\right)^2$$
for some positive integer $N_1$.
 Clearly,  $m$ divides
$N_1$.

 Clearly, $D$ contains a (maximal) field extension $L/E'$ of
degree $\frac{N_1}{m}$ and $\End^0(X,i)\cong \M_m(D)$ contains
every field extension $T/L$ of degree $m$. This implies that
$$\End^0(X) \supset\End^0(X,i)\supset T$$
and the number field $T$ has degree $[T:\Q]=[E':\Q]\cdot
\frac{N_1}{m} \cdot m=[E:\Q]e N_1$. But $[T:\Q]$ must divide
$2\dim(X)$ (see \cite[proposition 2 on p. 36]{Shimura}); if the
equality holds then $X$ is an abelian variety of CM-type. This
implies that $e N_1$ divides $d=\frac{2\dim(X)}{[E:\Q]}$. It
follows that $(eN_1)^2$ divides $d^2$; if the equality holds then
$[T:\Q]=2\dim(X)$ and therefore $X$ is an abelian variety of
CM-type. But $(eN_1)^2=e^2 N_1^2=e (e
N_1^2)=eN=e\cdot\dim_E(\End^0(X,i))$. This implies that
$\dim_E(\End^0(X,i)) \le \frac{d^2}{e} \le d^2$, which proves (i).

Assume now that
 $\dim_E(\End^0(X,i))=d^2$.  Then $e=1$ and
$$(eN_1)^2=r^2, N_1=d,\ [T:\Q]=[E:\Q]e N_1=[E:\Q] d=2\dim(X);$$
in particular, $X$ is an abelian variety of CM-type. In addition,
since $e=1$, we have $E'=E$, i.e. $\End^0(X,i)$ is a {\sl central}
simple $E$-algebra. We also have $C\subset E$ and
$$\dim_{E}(D)=\dim_{E'}(D)=\left(\frac{N_1}{m}\right)^2=\left(\frac{d}{m}\right)^2.$$
Since $E$ is the center of $D$, it is also  the center of the
matrix algebra $\M_m(D)$. Clearly, there exist an abelian variety
$Z$ over $K_a$, an embedding $j:D \hookrightarrow \End^0(Z)$ and
an isogeny $\psi:Z^m \to X$ such that the induced isomorphism
$$\psi_{*}: \End^0(Z^m) \cong \End^0(X),\ u\mapsto \psi u \psi^{-1}$$ maps
$j(\M_m(D)):=\M_m(j(D)) \subset \M_m(\End^0(Z))=\End^0(Z^m)$ onto
$\End^0(X,i)$. Since $E$ is the center of $\M_m(D)$ and $i(E)$ is
the center of $\End^0(X,i)$, the isomorphism $\psi_{*}$ maps $j(E)
\subset j(\M_m(D))=\M_m(j(D))\subset \End^0(Z^m)$
 onto $i(E) \subset \End^0(X)$. In other words,
 $\psi_{*}j(E)=i(E)$.
 It follows that there exists an
 automorphism $\sigma$ of the
 field $E$ such that
 $i=\psi_{*}j\sigma$ on $E$. This implies that if we put
 $k:=j\sigma:E \hookrightarrow \End^0(Z)$
 then $\psi \in\Hom((Z^m,k^{(m)}), (X,\psi))$.

 Clearly, $k(E)=j(E)$ and therefore $j(D)\subset \End^0(Z,k)$.
 Since $\M_m(\End^0(Z,k))\cong \End^0(X,i)\cong \M_m(D)$, the
 dimension arguments imply that $j(D)= \End^0(Z,k)$ and therefore
 $\End^0(Z,k)\cong D$ is a division algebra. We have
 $$\dim(Z)=\frac{\dim(X)}{m}, \quad
 \dim_E(D)=\left(\frac{d}{m}\right)^2=\left(\frac{2\dim(X)}{[E:\Q]m}\right)^2=\left(\frac{2\dim(Z)}{[E:\Q]}\right)^2.$$

Let $B$ be an absolutely {\sl simple} abelian variety over $K_a$
such that $X$ is isogenous to a self-product $B^r$ of $B$ where
the positive integer $r=\frac{\dim(X)}{\dim(B)}$. Then $\End^0(B)$
is a central division algebra over $C$; we define a positive
integer $g_0$ by $\dim_C(\End^0(B))=g_0^2$. Since $\End^0(X)$
contains a field of degree $2\dim(X)$, it follows from
Propositions 3 and 4 on pp. 36--37 in \cite{Shimura} (applied to
$A=X,K=C,g=g_0,m=\dim(B),f=[C:\Q]$) that $2\dim(B)=[C:\Q]\cdot
g_0$. Let $T_0$ be a maximal subfield in the $g_0^2$-dimensional
central division algebra $\End^0(B)$. Well-known properties of
maximal subfields of division algebras imply that $T_0$ contains
the center $C$ and $[T_0:C]=g_0$. It follows that $[T_0:\Q]=[C:\Q]
[T_0:C]=[C:\Q] \cdot g_0=2\dim(B)$ and therefore $\End^0(B)$
contains a field of degree $2\dim(B)$. This implies that $B$ is an
absolutely simple abelian variety of CM-type; in terminology of
\cite{Oort2}, $B$ is an absolutely simple abelian variety with
{\sl sufficiently many complex multiplications}.

 Assume now that $\fchar(K_a)=0$. We need to check that $2\dim(Z)=[E:\Q]$ and $E$ contains a
 CM-field.
 Indeed, since $D$ is a division algebra, it follows from Albert's classification
 \cite{MumfordAV,Oort} that $\dim_{\Q}(D)$
divides $2\dim(Z)=\frac{2\dim(X)}{m}=[E:\Q]\frac{d}{m}$. On the
other hand,
$\dim_{\Q}(D)=[E:\Q]\dim_{E}(D)=[E:\Q]\left(\frac{d}{m}\right)^2$.
Since  $m$ divides $d$, we conclude that $\frac{d}{m}=1$, i.e.,
$\dim_{E}(D)=1,  D=E,  2\dim(Z)=[E:\Q]$. In other words,
$\End^0(Z)$ contains the field $E$ of degree $2\dim(Z)$. It
follows from Theorem 1 on p. 40 in \cite{Shimura} (applied to
$F=E$) that $E$ contains a CM-field.

Now let us drop the assumption about $\fchar(K_a)$ and assume
instead that $E$  does {\sl not} contain a CM subfield. It follows
that
 $\fchar(K)>0$.  Since $C$ lies in $E$, it is totally real.
Since $B$ is an absolutely simple abelian variety with {\sl
sufficiently many complex multiplications}  it is isogenous to an
absolutely simple abelian variety $W$ defined over a finite field
\cite{Oort2} and $\End^0(B)\cong \End^0(W)$. In particular, the
center of $\End^0(W)$ is isomorphic to $C$ and therefore is a
totally real number field.  It follows from the Honda--Tate theory
\cite{Tate2} that $W$ is a supersingular elliptic curve and
therefore  $B$ is also a supersingular elliptic curve. Since $X$
is isogenous to $B^r$, it is  a supersingular abelian variety.

Now let us consider the case of arbitrary $X$. Applying the
already proven case of  Theorem \ref{maxE}(i) to each $X_s$, we
conclude that
$$\dim_E(\End^0(X_s,i)) \le
\left(\frac{2\dim(X_s)}{[E:\Q]}\right)^2.$$  Applying (5), we
conclude that
\begin{gather*}
\dim_E(\End^0(X,i)) =\sum_{s\in\I}\dim_E(\End^0(X_s,i_s))\le\\
\sum_{s\in\I}\left(\frac{2\dim(X_s)}{[E:\Q]}\right)^2\le
\frac{(2\sum_{s\in\I}
\dim(X_s))^2}{[E:\Q]^2}=\frac{(2\dim(X))^2}{[E:\Q]^2}.
\end{gather*}
 It follows
that if the equality
$\dim_E(\End^0(X,i))=\frac{(2\dim(X))^2}{[E:\Q]^2}$ holds then the
set $\I$ of indices $s$ is a singleton, i.e. $X=X_s$ is isogenous
to a self-product of  an absolutely simple abelian variety.

\end{proof}

\end{document}